\documentclass{gtart}
\input gtmonout
\volumenumber{1}
\volumeyear{1998}
\volumename{The Epstein birthday schrift}
\usepackage{rlepsf}
\pagenumbers{419}{450}
\received{13 November 1997}
\revised{15 October 1998}
\published{25 October 1998}
\papernumber{21}

\newtheorem{lem}{Lemma}[section]
\newtheorem{thm}[lem]{Theorem}
\newtheorem{cor}[lem]{Corollary}
\newtheorem{prop}[lem]{Proposition}
\theoremstyle{definition}
\newtheorem{define}[lem]{Definition}
\newtheorem{rem}{Remark}
\begin{document}

\title{Divergent sequences of Kleinian groups}
\author{Ken'ichi Ohshika}
\address{Graduate School of Mathematical Sciences\\
University of Tokyo\\
Komaba, Meguro-ku, Tokyo 153, Japan}
\email{ohshika@ms.u-tokyo.ac.jp}

\begin{abstract}
One of the basic problems in studying topological structures of deformation
spaces for Kleinian groups is to find a criterion to distinguish convergent
sequences from divergent sequences.
In this paper, we shall give a sufficient condition for sequences of
Kleinian groups isomorphic to surface groups to diverge in the deformation
spaces.
\end{abstract}

\primaryclass{57M50}\secondaryclass{30F40}

\keywords{Kleinian group, hyperbolic $3$--manifold,
deformation space}
\asciikeywords{Kleinian group, hyperbolic 3-manifold,
deformation space}

\maketitle

\cl{\em\small Dedicated to Prof David Epstein on the occasion of  
his 60th birthday}

\section{Introduction}

The deformation space of a Kleinian group $\Gamma$ is the space of faithful
discrete representations of $\Gamma$ into $PSL_2{\bf C}$ preserving
parabolicity modulo conjugacy.
It is one of the important aspects of Kleinian group theory  to study the
structures of deformation spaces.
The first thing that was studied among the structures of deformation spaces
was that of subspaces called quasi-conformal deformation spaces.
By  works of Ahlfors, Bers, Kra, Marden and Sullivan among others, the
topological types and the parametrization of quasi-conformal deformation
spaces are completely determined  using the theory of quasi-conformal
mappings and the ergodic theory on the sphere
(\cite{Ah*}, \cite{Be}, \cite{Ma}, \cite{Su0}).
On the other hand, the total deformation spaces are less understood.
A recent work of Minsky \cite{Mi} makes it possible to determine the
topological structure of  the total deformation space  completely in the
case of once-puncture torus groups.
The other cases are far from complete understanding.
Although very rough topological structures, for instance the connected
components of deformation spaces can be understood by virtue of recent
works of Anderson--Canary and Anderson--Canary--McCullough, more detailed
structures like the frontier of quasi-conformal deformation spaces are not
yet known even in the case of surface groups with genus greater than $1$.

A first step to understand the topological structure of the deformation
space of a Kleinian group $\Gamma$ is to give a criterion for a sequence
$\{\Gamma_i\}$ in the deformation space to converge or diverge.
In this paper, we shall consider the  simplest case when the group $\Gamma$
is isomorphic to a hyperbolic surface group $\pi_1(S)$ and has no
accidental parabolic elements.
In this case, $\Gamma_i$ is either quasi-Fuchsian or a totally degenerate
b--group, or a totally doubly degenerate group.
Hence by taking a subsequence, we have only to consider the following three
cases: all of the $\{\Gamma_i\}$ are quasi-Fuchsian, or totally degenerate
b--groups, or totally doubly degenerate groups.
For such groups, some conditions for sequences to converge are given for
example in Bers~\cite{Be}, Thurston~\cite{Th2} and Ohshika~\cite{Oh2}.
Thurston's convergence theorem is called the double limit theorem. 
The purpose of this paper is to give a sufficient condition for sequences
to diverge in the deformation space, which is in some sense complementary
to the condition of the double limit theorem.

Before explaining the content of our main theorem, let us recall that a
Kleinian group  isomorphic to a hyperbolic surface group without accidental
parabolic elements has two pieces of information describing the structures
near ends as follows.
When such a Kleinian group $\Gamma$ is quasi-Fuchsian, by the Ahlfors--Bers
theory, we get a pair of points in the Teichmu\"{u}ller space ${\cal T}(S)$
corresponding to the group.
In the case when $\Gamma$ is a totally degenerate b--group, as there is one
end of the non-cuspidal part $({\bf H}^3/\Gamma)_0$ which is geometrically
finite, we have a point in the Teichm\"{u}ller space.
In addition, the geometrically infinite end of $({\bf H}^3/\Gamma)_0$
determines an ending lamination which is defined uniquely up to changes of 
 transverse measures.
Finally in the case when $\Gamma$ is a totally doubly degenerate group,
$({\bf H}^3/\Gamma)_0$ have two geometrically infinite ends, and we have a
pair of measured laminations  which are ending laminations of  the two 
ends.
We shall define an end invariant of such a group $\Gamma$ to be a pair
$(\chi, \upsilon)$ where each factor is either a point of the
Teichm\"{u}ller space or a projective lamination represented by an ending
lamination, which gives the information on one of the ends.

The statement of our main theorem is as follows.
Suppose that we are given  a sequence of Kleinian groups
$(\Gamma_i,\phi_i)$ in the parabolicity-preserving deformation space
$AH_p(S)$ of Kleinian groups isomorphic to $\pi_1(S)$ for a hyperbolic
surface $S$.
Suppose moreover that the end invariants $(\chi_i , \upsilon_i)$ have the
following property:
Either in the Thurston compactification or in the projective lamination
space, $\{\chi_i\}$ and $\{\upsilon_i\}$ converge to maximal and connected
projective laminations with the same support.
Then the sequence $\{(\Gamma_i, \phi_i)\}$ does not converge in $AH_p(S)$.

To understand the meaning of this theorem, let us contrast it with
Thurston's double limit theorem.
For simplicity, we only consider the case when $\Gamma_i$ is a
quasi-Fuchsian group for the time being.
By Ahlfors--Bers theory, a sequence of quasi-Fuchsian groups $\{(\Gamma_i,
\phi_i)\}$ corresponds to a sequence of pairs of marked hyperbolic
structures $\{(m_i, n_i)\}$ on $S$.
Consider the case when both $m_i$ and $n_i$ diverge in the Teichm\"{u}ller
space and their limits in the Thurston compactification are projective
laminations $[\mu]$ and $[\nu]$ respectively. 
The double limit theorem asserts that if $\mu$ and $\nu$ fill up $S$, viz.,
any measured lamination has non-zero intersection number with either $\mu$
or $\nu$, then the sequence $\{(\Gamma_i, \phi_i)\}$ converges in the
deformation space passing through a subsequence if necessary.
The situation of our theorem is at the opposite pole to that of the double
limit theorem.
We assume in our theorem that $\mu$ and $\nu$ are equal except for the
transverse measures and that they are maximal and connected.

We can see the assumption of maximality is essential by taking look at an
example of Anderson--Canary \cite{AC}.
They constructed an example of  quasi-Fuchsian groups converging in
$AH_p(S)$ which correspond to  pairs of marked hyperbolic structures
$(m_i,n_i)$ such that $\{m_i\}$ and $\{n_i\}$ converge to the same point in
${\cal PL}(S)$.
In this example, the support of the limit projective lamination is a simple
closed curve, far from being maximal.

The proof of our theorem is based on an argument sketched in Thurston
\cite{Th} which was used to prove his theorem stating  that sequences of
Kleinian groups isomorphic to surface groups which converge algebraically
to Kleinian groups without accidental parabolic elements converge strongly.
We shall give a detailed proof of this theorem in the last section as an
application of our theorem.

The original version of this paper was written during the author's stay in
University of Warwick for the symposium ``Analytic and geometric aspects of
hyperbolic spaces''.
The author would like to express his gratitude to the organizers of the
symposium, Professors David Epstein and Caroline Series for inviting him
there and giving him a lot of mathematical stimuli.  

\section{Preliminaries}
\label{prem}

Kleinian groups are discrete subgroups of the Lie group $PSL_2{\bf C}$
which is the group of conformal automorphisms of the $2$--sphere $S^2$ and
the orientation preserving isometry group of the hyperbolic 3--space ${\bf
H}^3$.
A Kleinian group acts conformally on  $S^2$ and discontinuously on ${\bf
H}^3$ by isometries.
In this paper, we always assume that Kleinian groups are torsion free.
For a torsion-free Kleinian group $\Gamma$, the quotient ${\bf H}^3/\Gamma$
is a complete hyperbolic 3--manifold.

Let $\Gamma$ be a Kleinian group, which is regarded as  acting  on  $S^2$.
The subset of $S^2$ which is the closure of the set consisting of the fixed
points of non-trivial elements in $\Gamma$, is called the limit set of
$\Gamma$, and denoted by $\Lambda_\Gamma$.
The limit set $\Lambda_\Gamma$ is invariant under the action of $\Gamma$.
The complement of $\Lambda_\Gamma$ is called the region of discontinuity of
$\Gamma$ and denoted by $\Omega_\Gamma$.
The group $\Gamma$ acts on $\Omega_\Gamma$ properly discontinuously.
If $\Gamma$ is finitely generated, the quotient $\Omega_\Gamma/\Gamma$ is a
Riemann surface of finite type (ie a disjoint union of finitely many
connected Riemann surfaces of finite genus with finitely many punctures) by
Ahlfors' finiteness theorem \cite{Ah}.

A homeomorphism $\omega \co  S^2 \rightarrow S^2$ is said to be
quasi-conformal if it has an $L^2$--distributional derivative and there
exists a function $\mu \co  S^2 \rightarrow {\bf C}$ called a Beltrami
coefficient whose essential norm is strictly less than $1$, such that
$\omega_{\overline{z}}= \mu \omega_{z}$.
If the Beltrami coefficient $\mu$ for $\omega$ satisfies the condition $\mu
\circ \gamma(z) \overline{\gamma}'(z) / \gamma'(z) = \mu(z)$ for every
$\gamma \in \Gamma$, then the conjugate $\omega \Gamma \omega^{-1}$ is
again a Kleinian group.
A Kleinian group obtained by such a fashion from $\Gamma$ is called a
quasi-conformal deformation of $\Gamma$.
By identifying two quasi-conformal deformations which are conformally
conjugate, and giving the topology induced from the representation space,
we obtain the quasi-conformal deformation space of $\Gamma$, which we shall
denote by $QH(\Gamma)$.
A quasi-conformal deformation of $\Omega_\Gamma/\Gamma$ can be extended to
that of $\Gamma$. This gives rise to a continuous map $\rho\co {\cal
T}(\Omega_\Gamma/\Gamma) \rightarrow QH(\Gamma)$.
By the works of Ahlfors, Bers, Kra, Marden and Sullivan among others, it is
known that when $\Gamma$ is finitely generated, $\rho$ is a covering map,
and that especially if $\Gamma$ is isomorphic to a surface group (or more
generally if $\Gamma$ satisfies the condition ($*$) introduced by 
Bonahon~\cite{Bo}), then $\rho$ is a homeomorphism. The inverse of 
$\rho$ is denoted by $Q$.

Let $\Gamma$ be a finitely generated Kleinian group.
We shall define the deformation space of $\Gamma$.
An element $\gamma$ of $PSL_2{\bf C}$ is said to be parabolic if it is
conjugate to a parabolic element $\left( \matrix{1& 1 \hfill\cr 0&1\hfill\cr}
\right)$.
The deformation space of $\Gamma$, denoted by $AH_p(\Gamma)$, is the space
of faithful discrete representations of $\Gamma$ into $PSL_2{\bf C}$
preserving the parabolicity modulo conjugacy with the quotient topology
induced from the representation space.
We shall often denote an element (ie an equivalence class of groups) in
$AH_p(\Gamma)$ in a form $(G,\phi)$ where $\phi$ is a faithful discrete
representation with the image $G$ which represents the equivalence class.
The quasi-conformal deformation space $QH(\Gamma)$ is regarded as a
subspace of $AH_p(\Gamma)$.

Let $C(\Lambda_\Gamma)$ be the intersection of ${\bf H}^3$ and the convex
hull of the limit set $\Lambda_\Gamma$ in the Poincar\'{e} ball  ${\bf H}^3
\cup S^2_\infty$.
As $C(\Lambda_\Gamma)$ is $\Gamma$--invariant, $C(\Lambda_\Gamma)$ can be
taken quotient by $\Gamma$ and gives rise to a closed convex set
$C(\Lambda_\Gamma)/\Gamma$ in ${\bf H}^3/\Gamma$, which is called the
convex core of ${\bf H}^3/\Gamma$.
The convex core is the minimal closed convex set of ${\bf H}^3/\Gamma$
which is a deformation retract.
A  Kleinian group $\Gamma$ is said to be geometrically finite if it is
finitely generated and if the convex core of ${\bf H}^3/\Gamma$ has finite
volume, otherwise it is geometrically infinite.
When $\Gamma$ is geometrically finite, $QH(\Gamma)$ is an open subset of
$AH_p(\Gamma)$.

For a sequence $\{\Gamma_i\}$ of Kleinian groups, its geometric limit  is
defined as 
follows.
\begin{define}
\label{geometric limit}
A Kleinian group $H$ is called the geometric limit of $\{\Gamma_i\}$ if
every element 
of $H$ is the limit of a sequence $\{\gamma_i\}$ for $\gamma_i \in
\Gamma_i$, and 
the limit of any convergent sequence $\{\gamma_{i_j} \in \Gamma_{i_j}\}$
for a subsequence $\{\Gamma_{i_j}\} \subset \{\Gamma_i\}$ is contained in
$H$.
\end{define}

The geometric limit of non-elementary Kleinian groups is also a Kleinian group.
We call a limit in the deformation space an algebraic limit to distinguish
it from a geometric limit.
We also call the first factor of a limit in the deformation space, ie
the Kleinian group which is the image of the limit representation, an
algebraic limit.
Suppose that $\{(\Gamma_i, \phi_i)\}$ converges in $AH_p(\Gamma)$ to
$(\Gamma', \phi)$.
Then there is a subsequence of $\{\Gamma_i\}$ converging to a Kleinian
group $H$ geometrically.
Moreover, the algebraic limit $\Gamma'$ is contained in the geometric limit $H$.
(Refer to J{\o}rgensen--Marden \cite{JM} for the proofs of these facts.)
When the algebraic limit $\Gamma'$ coincides with the geometric limit $H$,
we say that the sequence $\{(\Gamma_i,\phi_i)\}$ converges to
$(\Gamma',\phi)$ strongly.

When $\{\Gamma_i\}$ converges geometrically to $H$, there exists a framed
$(K_i,r_i)$--app\-rox\-imate isometry defined below between ${\bf H}^3/\Gamma_i$
and ${\bf H}^3/H$ with base-frames which are the projections of a
base-frame on  a point in ${\bf H}^3$ where $K_i \rightarrow 1$ and $r_i
\rightarrow \infty$ as $i\rightarrow \infty$.
(See Canary--Epstein--Green~\cite{CEG}).)

\begin{define}
Let $(M_1,e_1)$ and $(M_2,e_2)$ be two Riemannian 3--manifolds with
base-frame whose base-frames are based at $x_1 \in M_1$, and $x_2 \in M_2$
respectively.
A $(K,r)$--approximate isometry between $(M_1,e_1)$ and $(M_2,e_2)$ is a 
diffeomorphism from  $(X_1,x_1)$ to  $(X_2,x_2)$ for subsets $X_1,X_2$ of
$M_1,M_2$ 
containing the $r$--balls centred at $x_1, x_2$ such that $df(e_1) =e_2$ and
$$d_{M_1}(x,y)/K \leq d_{M_2}(f(x),f(y)) \leq Kd_{M_1}(x,y)$$ for any $x,y
\in X_1$.
\label{approximate isometry}
\end{define}

Let $\{(M_i,v_i)\}$ be a sequence of hyperbolic 3--manifolds with base-frame.
We say that $(M_i,v_i)$ converges  geometrically (in the sense of Gromov)
to  a hyperbolic 3--manifold with base-frame $(N,w)$ when for any large $r$
and $K>1$ there exists an integer $i_0$ such that there exists a
$(K,r)$--approximate isometry between $(M_i,v_i)$ and $(N,w)$ for $i \geq
i_0$.
As described above, by choosing base-frames which are the images of a fixed
base-frame in ${\bf H}^3$, the sequence of ${\bf H}^3/\Gamma_i$ with the
base-frame converges geometrically to ${\bf H}^3/H$ with the base-frame
when $\Gamma_i$ converges to $H$ geometrically.

Let $M={\bf H}^3/\Gamma$ be a complete hyperbolic 3--manifold.
A parabolic element of $\Gamma$ is contained in a maximal parabolic
subgroup, which is  isomorphic to either ${\bf Z}$ or ${\bf Z} \times {\bf
Z}$ and corresponds to a cusp of $M$.
This is derived from Margulis' lemma.
By deleting mutually disjoint neighbourhoods of the cusps of $M$, we obtain
a  non-cuspidal part of $M$, which we shall denote by $M_0$.
We delete the cusp neighbourhoods where the injectivity radius is less than
$\epsilon$ for some universal constant $\epsilon >0$ so that this procedure
of deleting cusp neighbourhoods is consistent among all the  hyperbolic
3--manifolds.
The non-cuspidal part $M_0$ is a 3--manifold whose boundary component is
either a torus or an open annulus.

By  theorems of Scott~\cite{Sc} and McCullough~\cite{Mc}, there exists a
submanifold $C(M)$ of $M_0$ such that $(C(M), C(M) \cap \partial M_0)$ is
relatively homotopy equivalent to $(M_0, \partial M_0)$ by the inclusion,
which is called a core of $M$.
An end of $M_0$ is said to be geometrically finite if some neighbourhood of
the end contains no closed geodesics, otherwise it is called geometrically
infinite.
A geometrically infinite end $e$ is called geometrically infinite tame (or
simply degenerate) if that end faces an incompressible frontier component
$S$ of a core and there exists a sequence of simple closed curves
$\{\gamma_i\}$ on $S$ such that the closed geodesic in $M$ homotopic to
$\gamma_i$ tends to the end $e$ as $i \rightarrow \infty$.
(In this paper we use this term only when every component of the frontier
of the core is incompressible.)
A Kleinian group $\Gamma$ is geometrically finite if and only if every end
of $({\bf H}^3/\Gamma)_0$ is geometrically finite.

In this paper, we shall consider sequences of Kleinian groups isomorphic to
surface groups.
Let $S$ be a hyperbolic surface of finite area.
We call punctures of $S$ cusps.
We denote by $AH_p(S)$ the space of Kleinian groups modulo conjugacy which
are isomorphic to $\pi_1(S)$ by isomorphisms mapping elements represented
by cusps to parabolic elements.
We can also identify this space $AH_p(S)$ with the deformation space of a
Fuchsian group $G$ such that ${\bf H}^2/G = S$.
Let $(\Gamma , \phi)$ be a class in $AH_p(S)$.
We say that a parabolic element $\gamma \in \Gamma$ is accidental parabolic
when $\phi^{-1}(\gamma)$ does not correspond to  a cusp of $S$.
Assume that  $(\Gamma, \phi)$ in $AH_p(S)$ has no accidental parabolic element.
Then the non-cuspidal part $({\bf H}^3/\Gamma)_0$ has only two ends since
one can see that a core is homeomorphic to $S \times I$ and has exactly two
frontier components.
Therefore in this case, $\Gamma$ is either (1) a quasi-Fuchsian group,
ie geometrically finite and the limit set $\Lambda_\Gamma$ is
homeomorphic to the circle or (2) a totally degenerate b--group, ie
$\Omega_\Gamma$ is connected and simply connected, and $({\bf
H}^3/\Gamma)_0$ has one geometrically finite end and one geometrically
infinite end, or (3) a totally doubly degenerate group, ie
$\Omega_\Gamma = \emptyset$, and $({\bf H}^3/\Gamma)_0$ has two
geometrically infinite tame ends.
Recall that a Kleinian group is called a b--group when its region of
discontinuity has a unique invariant component, which is simply connected.

For a hyperbolic surface $S= {\bf H}^2/\Gamma$, we denote the
quasi-conformal deformation space of $\Gamma$ by $QF(S)$.
This space consists of quasi-Fuchsian groups isomorphic to $\pi_1(S)$ by
isomorphisms taking elements representing cusps to parabolic elements.
By the Ahlfors--Bers theory, there is a homeomorphism $Q\co  QF(S) \rightarrow
{\cal T}(S) \times {\cal T}(\overline{S})$, which we shall call the
Ahlfors--Bers homeomorphism.
Here $\cal{T}(\overline{S})$ denotes the Teichm\"{u}ller space of the
``complex conjugate" of $S$.
This can be interpreted as the space of marked hyperbolic structures on $S$
such that the complex conjugate of the corresponding complex structure is
equal to the structure on the second component of $\Omega_\Gamma/\Gamma$.
We identify $\cal{T}(\overline{S})$ with $\cal{T}(S)$ by the above
correspondence from now on.
By this correspondence,  the Fuchsian representations of $\pi_1(S)$ are
mapped onto the diagonal of ${\cal T}(S) \times {\cal T}(S)$. 

Thurston introduced a natural compactification of a Teichm\"{u}ller space
in \cite{Th0}, which is called the Thurston compactification nowadays.
Let $S$ be a hyperbolic surface of finite area.
Let ${\cal S}$ denote the set of free homotopy classes of simple closed
curves on $S$.
Let $P {\bf R}_+^{{\cal S}}$ denote the projective space obtained from the
space ${\bf R}^{\cal S}_+$ of non-negative functions on ${\cal S}$.
We endow $P {\bf R}_+^{{\cal S}}$ with the quotient topology of the weak
topology on ${\bf R}^{{\cal S}}_+ \setminus \{0\}$.
The Teichm\"{u}ller space $ {\cal T}(S)$ is embedded in $P{\bf R}_+^{{\cal
S}}$ by taking $g \in \cal{T}(S)$ to the class represented by a function 
whose value at $s \in {\cal S}$ is the length of the closed geodesic in the
homotopy class.
The closure of the image of ${\cal T}(S)$ in $P{\bf R}_+^{{\cal S}}$ is
homeomorphic to the ball and defined to be the Thurston compactification of
${\cal T}(S)$.
The boundary of ${\cal T}(S)$ corresponds to ``the space of projective
laminations" in the following way.

A compact subset of $S$ consisting of disjoint simple geodesics is called a
geodesic lamination.
A geodesic lamination endowed with a transverse measure which is invariant
under a homotopy along leaves is called a measured lamination.
The subset of a measured lamination $\lambda$ consisting of the points $x
\in \lambda$ such that any arc containing $x$ at the interior has a
positive measured with respect to the transverse measure is called the
support of $\lambda$.
We can easily see that the support of a measured lamination $\lambda$ is a
geodesic lamination.
The set of measured laminations with the weak topology with respect to
measures on  finite unions  of arcs is called the measured lamination space
and denoted by ${\cal ML}(S)$.
The set of simple closed geodesics with positive weight is dense in ${\cal
ML}(S)$.
For a measured lamination $(\lambda , \mu)$, where $\mu$ denotes the
transverse measure, and a homotopy class of simple closed curves $\sigma$,
we  define their intersection number $i(\lambda, \sigma)$ to be $\inf_{s
\in \sigma} \mu(s)$.
(We also use the notation $i(\lambda, s)$ to denote $i(\lambda, [s])$.)
By defining the value at $\sigma \in {\cal S}$ to be $i(\lambda, \sigma)$,
we can define a map $\iota \co   {\cal ML}(S) \rightarrow {\bf R}_+^{{\cal
S}}$.
By projectivising the both spaces, we have a map $\overline{\iota} \co  {\cal
PL}(S) \rightarrow P{\bf R}_+^{{\cal S}}$, where ${\cal PL}(S)$ denotes the
projectivization of ${\cal ML}(S)$, ie $({\cal ML}(S) \setminus 
\{\emptyset\})/(0,\infty)$.
It can be proved that in fact $\overline{\iota}$ is an embedding  and
coincides with the boundary of the image of ${\cal T}(S)$, that is, the
boundary of the Thurston compactification of ${\cal T}(S)$.
Refer to Fathi et al \cite{Fa} for further details of these facts.

Let $e$ be a geometrically infinite tame end of the non-cuspidal part of a
hyperbolic 3--manifold $M$, which faces a frontier component $\Sigma$ of a
core.
 From now on, we always assume that every frontier component of a core is
incompressible in $M$.
By the definition of geometrically infinite tame end, there exists a
sequence of simple closed curves $\{\gamma_i\}$ on $\Sigma$ such that the
closed geodesic homotopic to $\gamma_i$ tends to $e$ as $i \rightarrow
\infty$.
Consider the sequence $\{[\gamma_i]\}$ (the projective classes represented
by $\{\gamma_i\}$) in ${\cal PL}(\Sigma)$.
(We identify $\gamma_i$ with the closed geodesic homotopic to $\gamma_i$
with respect to some fixed hyperbolic structure on $\Sigma$.)
Since ${\cal PL}(\Sigma)$ is compact, the sequence $\{[\gamma_i]\}$
converges to a projective lamination $[\lambda] \in {\cal PL}(\Sigma)$
after taking a subsequence.
Such a measured lamination $\lambda$ is called an ending lamination of $e$.
(The original definition is due to Thurston~\cite{Th}.)
An ending lamination is maximal (ie it is not a proper sublamination of
another measured lamination), and connected.
(Thurston~\cite{Th}, see also Ohshika~\cite{Oh1}.)
If both $\lambda$ and $\lambda'$ are ending laminations of an end $e$,
their intersection number $i(\lambda,\lambda')$ is equal to $0$
(essentially due to Thurston~\cite{Th} and Bonahon~\cite{Bo}).
We shall give a proof of this fact, based on Bonahon's result in section~3.
By the maximality, this implies that $|\lambda|= |\lambda'|$ where
$|\lambda|$ denotes the support of $\lambda$.

In this paper, we shall deal with a hyperbolic 3--manifold $M={\bf
H}^3/\Gamma$ with a homotopy equivalence $\phi \co  S \rightarrow M$
preserving cusps.
In this case, $M$ has a core which is homeomorphic to $S \times I$.
For a homotopy equivalence $\phi \co  S \rightarrow M$ and a lamination
$\lambda$, its image $\phi(\lambda)$ is homotopic to a unique  lamination
on $S\times \{t\}$  for both $t=0,1$.
When the measured lamination  homotopic to $\phi(\lambda)$ is an ending
lamination, we  say that $\phi(\lambda)$ represents an ending lamination.
For an end $e$ of $M$, the {\em end invariant} of $e$ is defined to be a
projective lamination $[\lambda]$ on $S$ such that $\phi(\lambda)$
represent an ending lamination of $e$ when $e$ is geometrically infinite,
and the point in the Teichm\"{u}ller space corresponding to the conformal
structure of the component of $\Omega_ \Gamma/\Gamma$ when $e$ is
geometrically finite.

Now, let $e_1, e_2$ be the two ends of $M_0$  which are contained in the
``upper complement" and the ``lower complement"  of a core respectively
with respect to the orientation give on $M$ and $S$.
We define the {\em end invariant} of $M={\bf H}^3/\Gamma$ to be a pair
$(\chi, \upsilon)$, where $\chi$ is the end invariant of $e_1$ and
$\upsilon$ that of $e_2$.
This means in particular that when $\Gamma$ is a quasi--Fuchsian group, the
end invariant is equal to $Q(\Gamma, \phi) \in {\cal T}(S) \times {\cal
T}(S)$, where $Q \co  QF(S) \rightarrow {\cal T}(S) \times {\cal T}(S)$ is the
Ahlfors--Bers map with the second factor ${\cal T}(\overline{S})$ identified
with ${\cal T}(S)$.

Let $S$ be a hyperbolic surface of finite area and $M$  a complete
hyperbolic 3--manifold.
A pleated surface $f\co  S\rightarrow M$ is a continuous map which is totally
geodesic in $S-\ell$ for some geodesic lamination $\ell$ on $S$ such that
the path metric induced by $f$ coincides with the hyperbolic metric on
$S$.
We say that a sequence of pleated surfaces with base point $\{f_i \co 
(S_i,x_i) \rightarrow (M_i,y_i)\}$ converges geometrically to  a pleated
surface with base point $f \co  (S,x) \rightarrow (M,y)$ when there are
$(K_i,r_i)$--approximate isometries $\rho_i$ between $(M_i, v_i)$ and
$(M,v)$, and $\overline{\rho}_i$ between $(S_i,w_i)$ and $(S,w)$ such that
$K_i \rightarrow 1$ and $r_i \rightarrow \infty$ as $i \rightarrow \infty$
and $\{
\rho_i \circ f_i \circ \overline{\rho}_i^{-1}\}$ converges to $f$ uniformly on
every compact subset of $S$, where $v_i, v, w_i,w$ are base-frames on $x_i,
x, y_i, y$ respectively.
The space of pleated surfaces has the following compactness property due to
Thurston whose proof can be found in Canary--Epstein--Green~\cite{CEG}.

\begin{prop}
\label{compact}
For any sequence of pleated surfaces with base point $\{f_i \co$\break
 $(S_i,x_i)
\rightarrow (M_i,y_i)\}$ such that the injectivity radius at $y_i$ is
bounded away from $0$ as $i \rightarrow \infty$, there exists a subsequence
which converges geometrically.
\end{prop}

We say that a (measured or unmeasured) geodesic lamination $\lambda $ on
$S$ is realized by a pleated surface $f$ when $\lambda$ is mapped totally
geodesically by $f$.
A measured lamination $\lambda$ lying  on a component of the frontier of a
core of $M$ represents an ending lamination of an end of $M_0$ if and only
if there is no pleated surface (homotopic to the inclusion) realizing
$\lambda$.
(This follows from Proposition~5.1 in Bonahon~\cite{Bo} which we shall cite
below as Proposition~\ref{5.1}.)

We shall use the following two results of Bonahon~\cite{Bo} several times
in this paper.
The first is Proposition~3.4 in his paper.

\begin{lem}[Bonahon]
\label{intersection}
Let $M$ be a complete hyperbolic 3--manifold.
Let $S$ be a properly embedded incompressible surface in the non-cuspidal
part $M_0$.
Then there exists a constant $C$ with the following property.
Let $\alpha^*, \beta^*$ be closed geodesics in $M$ which are homotopic to
closed curves $\alpha, \beta$ on $S$ by homotopies coming to the same side
of $S$, and are located at distance at least $D$ from $S$.
Suppose that neither $\alpha^*$ nor $\beta^*$ intersects
a Margulis tube whose axis is not itself, $\alpha^*$ or $\beta^*$.
Then we have $$i(\alpha, \beta) \leq Ce^{-D}{\rm length}(\alpha){\rm
length}(\beta) +2.$$
\end{lem}

The second is Proposition~5.1 in Bonahon's paper.
Before stating the proposition, we need to define some terms used there.
A train track on a surface $S$ is a graph with $C^1$--structure such that
all edges coming to a vertex are tangent mutually there.
Furthermore we impose the condition that there is no component of the
complement which  is the interior of a monogon or a bigon or an annulus
without angle.
We call edges of a train track branches and vertices switches.
A regular neighbourhood of a train track $\tau$ can be foliated by arcs
transverse to $\tau$.
Such a neighbourhood is called a tied neighbourhood of $\tau$, and the arcs
are called ties.
We say that a geodesic lamination $\lambda$ is carried by a train track
$\tau$ when  a tied neighbourhood of $\tau$ can be isotoped to contain
$\lambda$ so that each leaf of $\lambda$ should be transverse to the ties.

When $\lambda$ is a measured lamination and carried by a train track
$\tau$, the transverse measure induces a weight system on the branches  of
$\tau$, by defining the weight of a  branch to be the measure of ties
intersecting the branch.
We can easily prove that such a weight system is uniquely determined by
$\lambda$ and $\tau$.
Conversely a weight system $w$ on a train track $\tau$ satisfying the
switch condition that the sum of weights on incoming branches and the sum
of those on outgoing branches coincides at each switch, determines a unique
measured lamination such that the weight system which it induces on  $\tau$
is equal to $w$.
Refer to Penner--Harer \cite{PH} for more precise definitions and
explanations for these facts.

A continuous map $f$ from a surface $S$ to a hyperbolic manifold $M$ is
said to be adapted to a tied neighbourhood $N_\tau$ of a train track $\tau$
on $S$ when each branch of $\tau$ is mapped to a geodesic arc in $M$ and
each tie of $N_\tau$ is mapped to a point.
Consider a map $f$ adapted to a tied neighbourhood of a train track $\tau$.
For a weight system $w$ on $\tau$, we  define the length of $f(\tau, w)$ to
be  $\sum w_b {\rm length}(f(b))$, where the sum is taken over all the
branches of $\tau$, and $w_b$ denotes the weight on $b$ assigned by $w$.
For a measured lamination $\lambda$ carried by $\tau$, if it induces a
weight system $w$ on  $\tau$, we define the length of $f(\lambda)$ to be
the length of $f(\tau, w)$.

For two branches $b, b'$ meeting at a switch $\sigma$ from opposite
directions, the exterior angle $\theta(f(b,b'))$ between $b, b'$ with
respect to $f$ is the exterior angle formed by $f(b)$ and $f(b')$ at
$f(\sigma)$.
The weight system $w$ determines the weight flowing from $b$ to $b'$.
Let $b_1, \dots , b_p$ and $b_1',\dots ,b_q'$ be the branches meeting at a
switch $\sigma$ with $b_1, \dots, b_p$ coming from one direction and $b_1',
\dots , b_q'$ from the other.
The exterior angle at $f(\sigma)$ is the sum of
$w_{k,l}\theta(f(b_k,b_l'))$ for all $k=1, \dots, p, l=1, \dots , q$, where
$w_{k,l}$ denotes the weight flowing from $b_k$ to $b'_l$.
The quadratic variation of angle at $f(\sigma)$ is the sum of
$w_{k,l}\theta^2f(b_k, b_l')$ in the same situation as above.
The total curvature of $f(\tau, w)$ is defined to be the sum of the
exterior angles at all the images of switches on $\tau$.
Similarly, the quadratic  variation of angle for $f(\tau,w)$ is defined to
be the sum of the quadratic variations of angle at all switches.

\begin{prop}[Bonahon]
\label{5.1}
Let $M$ be a complete hyperbolic 3--manifold and $S$ a hyperbolic surface of
finite type.
Let $\phi \co  S \rightarrow M$ be a continuous incompressible map taking
cusps to cusps, and $\lambda$ a measured lamination on $S$.
Then the one of the following two cases occurs and they are mutually exclusive.
\begin{enumerate}
\item For any $\epsilon>0$, there is a map $\phi_\epsilon$ homotopic
to $\phi$, which is adapted to a train track carrying $\lambda$ such
that ${\rm length}(\phi_\epsilon(\lambda))< \epsilon$.
\item For any $\epsilon$, there is a map $\phi_\epsilon$ homotopic to
$\phi$, which is adapted to a train track $\tau$ carrying $\lambda$ by a weight
system $\omega$, with the following property:
The total curvature and the quadratic variation of angle for
$\phi_\epsilon(\tau, w)$ are less than $\epsilon$.
Furthermore such a map $\phi_\epsilon$ satisfies the following:
There are $\delta >0, t<1$ such that $\delta \rightarrow 0, t \rightarrow
1$ as $\epsilon \rightarrow 0$, and for any simple closed curve $\gamma$
such that $[\gamma]$ is sufficiently close to $[\lambda]$ in ${\cal
PL}(S)$, the closed geodesic $\gamma^*$ homotopic to $\phi(\gamma)$ in $M$
has a part of length at least $t {\rm length}\phi_\epsilon(\gamma)$ which
lies within distance $\delta$ from $\phi_\epsilon(\gamma)$.
\end{enumerate}
\end{prop}

We can easily see that the first alternative exactly corresponds to the
case when $\lambda$ represents an ending lamination, and that  the second
alternative holds if and only if there is a pleated surface realizing
$\lambda$.
Taking this into account, the proposition implies in particular the following.
First, in the situation as in the proposition, $\phi(\lambda)$  represents
an ending lamination of an end of $M_0$ if and only if it is not realized
by a pleated surface homotopic to $\phi$ since the two alternatives are
exclusive.

Secondly, if $\lambda$ is an ending lamination, then any measured
lamination $\lambda'$ with the same support as $\lambda$ is also an ending
lamination.
This is because a train track carrying $\lambda$ also carries $\lambda'$
and if the condition (1) holds for $\lambda$, it equally holds for the
weight system corresponding to $\lambda'$.

There is another proposition which we shall make use of essentially in our
proof.
The proposition is an application of  Thurston's covering theorem which
originally appeared in \cite{Th} (see also \cite{Oh3} for its proof, and
Canary \cite{Ca} for its generalization).

\begin{prop}[Thurston]
\label{covering}
Let $S$ be a hyperbolic surface of finite area.
Let $\{(\Gamma_i, \phi_i)\}$ be a sequence of Kleinian groups in  $AH_p(S)$
converging to $(G,\psi)$.
Let $\Gamma_\infty$ be a geometric limit of $\{\Gamma_i\}$ after taking a
subsequence, and let $q \co  {\bf H}^3/G \rightarrow {\bf H}^3/\Gamma_\infty$
be the covering map associated with the inclusion $G \subset
\Gamma_\infty$.
Suppose that $({\bf H}^3/G)_0$ has a geometrically infinite end $e$.
Then there exists a neighbourhood $E$ of $e$  such that $q|E$ is a proper
embedding.
\end{prop}

\section{The main theorem}
Our main theorem on a sufficient condition for Kleinian groups isomorphic
to surface groups to diverge in the deformation spaces is the following.

\begin{thm}
\label{main}
Let $S$ be a hyperbolic surface of finite area.
Let $\{(\Gamma_i,\phi_i)\}$ be a sequence of Kleinian groups in $AH_p(S)$
with isomorphisms $\phi_i \co  \pi_1(S) \rightarrow \Gamma_i$ inducing 
homotopy equivalences $\phi_i \co  S \rightarrow {\bf H}^3/\Gamma_i$.
Let $(\chi_i, \upsilon_i)$ be an end invariant of $(\Gamma_i,\phi_i)$.
Suppose that $\{\chi_i\}$ and $\{\upsilon_i\}$ converge  in either  the
Thurston compactification of the Teichm\"{u}ller space  ${\cal T}(S)$ or
the projective lamination space ${\cal PL}(S)$ to maximal connected
projective laminations $[\mu]$ and $[\nu]$ with the same support.
Then $\{(\Gamma_i,\phi_i)\}$ does not converge  in $AH_p(S)$.
\end{thm}

Let us briefly sketch the outline of the proof of our main theorem.
Note that we can assume by taking a subsequence that all the $\Gamma_i$ are
the same type of the three;  quasi-Fuchsian groups or totally degenerate
b--groups or totally doubly degenerate groups.
We consider here only the case when all the $\Gamma_i$ are quasi-Fuchsian.
The proof is by reductio ad absurdum.
Suppose that our sequence $\{(\Gamma_i,\phi_i)\}$ converges in $AH_p(S)$.
Then we have the algebraic limit $(G,\psi)$ which is a subgroup of a
geometric limit $\Gamma_\infty$.
By applying the continuity of the length function on $AH_p(S) \times {\cal
ML}(S)$, which will be stated and proved in Lemma \ref{continuity}, we
shall show that $\psi(\mu)$ represents an ending lamination of an end
$e_\mu$ in $({\bf H}^3/G)_0$.
We shall  take a neighbourhood $E_\mu$ of $e_\mu$ which can be projected
homeomorphically by the covering map $q \co  {\bf H}^3/G \rightarrow {\bf
H}^3/\Gamma_\infty$ to a neighbourhood of an end of $({\bf
H}^3/\Gamma_\infty)_0$ using Proposition \ref{covering}. 
Let $S_0$ denote the non-cuspidal part of $S$.
We shall then show that deep inside $E_\mu$ there is an embedded surface
$f'(S_0)$ homotopic to $\psi|S_0$ such that every pleated surface homotopic
to $q \circ \psi$ touching $q \circ f'(S_0)$ is contained in $q(E_\mu)$.

By projecting $f'$ to ${\bf H}^3/\Gamma_\infty$ and pulling back by an
approximate isometry, we get  an embedded surface $f_i \co  S_0 \rightarrow
{\bf H}^3/\Gamma_i$ which is homotopic to $\phi_i$ converging to an
embedded surface $f_\infty \co  S_0 \rightarrow {\bf H}^3/\Gamma_\infty$
geometrically  which is the projection of $f'$.
By using a technique of interpolating pleated surfaces due to Thurston, we
shall show that there is a pleated surface $k_i \co  S \rightarrow {\bf
H}^3/\Gamma_i$ homotopic to $\phi_i$ which intersects $f_i(S_0)$ at an
essential simple closed curve.
These pleated surfaces converge geometrically to a pleated surface
$k_\infty \co  S' \rightarrow {\bf H}^3/\Gamma_\infty$, where $S'$ is an open
incompressible surface on $S$.
The condition that the limit surface $k_\infty$ touches $f_\infty(S)$
forces $k_\infty$ to be a pleated surface from  $S$, and to  be lifted to a
pleated surface to ${\bf H}^3/G$ which realizes a measured lamination with
the same support as $\mu$.
This will contradict the fact that $\psi(\mu)$ represents an ending lamination.

\section{Ending laminations and pleated surfaces}
\label{sec not b}
In this section, we shall prove lemmata basically due to Thurston which
will be used in the proof of our main theorem.

Throughout this section, $\{(\Gamma_i,\phi_i)\}$ denotes a sequence as in
Theorem~\ref{main}.
Suppose that $\{(\Gamma_i,\phi_i)\}$ converges to $(G,\psi)$ in $AH_p(S)$
where $\psi \co \pi_1(S) \rightarrow G$ is an isomorphism.
(Our proof of Theorem~\ref{main} is by reductio ad absurdum.
Therefore we assumed above the contrary of the conclusion of
Theorem~\ref{main}.)
We also use this symbol $\psi$ to denote the homotopy equivalence from $S$
to ${\bf H}^3/G$ corresponding to the isomorphism.
We can assume that $\phi_i$ converges to $\psi$ as representations by
taking conjugates if necessary.

Now let $\tilde{z} \in {\bf H}^3$ be a point and $\tilde{v}$ be a
frame based on $\tilde{z}$.  Then $\tilde{z}, \tilde{v}$ are projected
by the universal covering maps to $z_i,v_i$ of ${\bf H}^3/\Gamma_i$
and $z,v$ of ${\bf H}^3/G$.  Since we assumed that $\{\Gamma_i\}$
converges algebraically to $G$, we can assume by passing through a
subsequence that $\{\Gamma_i\}$ converges geometrically to a Kleinian
group $\Gamma_\infty$ which contains $G$ as a subgroup.  Let $v_\infty
, z_\infty$ be the images in ${\bf H}^3/\Gamma_\infty$ of $\tilde{v},
\tilde{z}$ by the universal covering map.

The hyperbolic manifolds with base frame $\{({\bf H}^3/\Gamma_i,v_i)\}$
converge in the sense of Gromov to $({\bf H}^3/\Gamma_\infty,v_\infty)$.
Let $ q\co  {\bf H}^3/G \rightarrow {\bf H}^3/\Gamma_\infty$ be the covering
associated with the inclusion $G \subset \Gamma_\infty$.
Then $q(z) = z_\infty$ and $dq(v) = v_\infty$.

Consider the case when at least one end $e$ of $({\bf H}^3/\Gamma_i)_0$ is
geometrically finite.
Let $\Sigma_i$ be the boundary components of the convex core of ${\bf
H}^3/\Gamma_i$ facing $e$ which corresponds to a component of the quotient
of the region of discontinuity $\Omega^0_{\Gamma_i}/\Gamma_i$.
Let $h_i \co  S\rightarrow \Sigma_i$  be  a homeomorphism homotopic to $\phi_i$.
Now by the assumption of Theorem \ref{main}, the marked conformal
structures of $\Omega^0_{\Gamma_i}/\Gamma_i$ converge to either $[\mu]$ or
$[\nu]$, say $[\mu]$.
Then we have the following.

\begin{lem}
\label{sequence}
There exist  an essential  simple closed curve $\gamma_i$ on $\Sigma_i$,
and a sequences of positive real numbers $\{r_i\}$ going to $0$ such that
$r_i {\rm length}_{\Sigma_i}(\gamma_i) \rightarrow 0$ and 
$\{r_i(h_i^{-1}(\gamma_i)) \in {\cal ML}(S)\}$  converges to  a measured
lamination with the same support as the measured lamination $\mu$, where we
regard $h_i^{-1}(\gamma_i)$  as an element in ${\cal ML}(S)$.
\end{lem}

\noindent
\proof
Let $m_i$ be the point in ${\cal T}(S)$ determined by the marked conformal
structure on $\Omega_{\Gamma_i}^0/\Gamma_i$.
By Sullivan's theorem proved in Epstein--Marden~\cite{EM}, the assumption in
Theorem~\ref{main} that $m_i \rightarrow [\mu]$  implies that the marked
hyperbolic structures $g_i$ on $S$ induced by $h_i$ from those on
$\Sigma_i$  as subsurfaces in ${\bf H}^3/\Gamma_i$ also converge to $[\mu]$
 as $i \rightarrow \infty$ in the Thurston compactification of ${\cal
T}(S)$.

Let $\gamma_i$ be the shortest essential closed curves on $\Sigma_i$ with
respect to the hyperbolic metrics induced from ${\bf H}^3/\Gamma_i$.
Consider the limit $[\mu_0]$ of $\{[h_i^{-1}(\gamma_i)]\}$  in ${\cal
PL}(S)$ passing through a subsequence if necessary.
Then there are bounded sequences of positive real numbers $r_i$ such that
$r_i h_i^{-1}(\gamma_i) \rightarrow \mu_0$  in ${\cal ML}(S)$.
Suppose that $i(\mu, \mu_0) \neq 0$.
Then by the ``fundamental lemma" 8-II-1 in Fathi--Laudenbach--Poenaru, we
should have ${\rm length}(r_ih_i^{-1}(\gamma_i)) \rightarrow \infty$.
On the other hand, since $\gamma_i$ is the shortest essential closed curve
with respect to $g_i$,\break  we see that ${\rm length}_{g_i}(h^{-1}(\gamma_i)) =
{\rm length}_{\Sigma_i} (\gamma_i)$ is bounded.
This implies that\break $r_i{\rm length}(h_i^{-1}(\gamma_i))$ is also bounded as
$ i \rightarrow \infty$, which is a contradiction.
Thus we have proved that $i(\mu, \mu_0)=0$.

As $\mu$ is assumed to be maximal and connected, this means that
$|\mu|=|\mu_0|$.
In particular  $\mu_0$  is not a simple closed curve, and we can see the
sequences $\{r_i\}$ must go to $0$ as $i \rightarrow \infty$.
\endproof

The next lemma, which asserts the continuity of the lengths of realized
measured laminations, appeared in Thurston~\cite{Th2}.
The following proof is based on Proposition \ref{5.1} due to Bonahon.
Soma previously suggested a possibility of such a proof.

\begin{lem}
\label{continuity}
Let $L \co  AH_p(S) \times {\cal ML}(S) \rightarrow {\bf R}$ be the function
such that\break $L((\Gamma, \phi),\lambda)$ is the length of the realization of
$\lambda$ on a pleated surface homotopic to $\phi$ when such a pleated
surface exists, otherwise set $L((\Gamma, \phi),\lambda) = 0$.
Then $L$ is continuous.
\end{lem}

\noindent
\proof
Let $\{(G_i,\psi_i)\} \in AH_p(S)$ be a sequence which converges to
$(G',\psi') \in AH_p(S)$, and let $\{\lambda_j\}$ be measured laminations
on $S$ converging to $\lambda'$.
We shall prove that $L$ is continuous at $((G',\psi'),\lambda')$.
We can take representatives for elements of the sequence so that  the
representations $\{\psi_i\}$ converge to $\psi'$.
Fix a base frame $\tilde{v}$ on ${\bf H}^3$ and let $w_i$ be the base frame
of ${\bf H}^3/G_i$ which is the projection of $\tilde{v}$ by the universal
covering map.
Since $G_i$ converges algebraically, the injectivity radius at the
basepoint under $w_i$ is bounded away from $0$ as $i \rightarrow \infty$.
By compactness of geometric topology (see Corollary 3.1.7 in
Canary--Epstein--Green~\cite{CEG}) and the diagonal argument, we can see that
for any large $r>0$ and small $\epsilon>0$, there exists $i_0$ such that
for any $i > i_0$, there exists a Kleinian group $H'$ containing $G'$ and a
$((1 + \epsilon),r)$--approximate isometry $\rho_i \co  B_r({\bf H}^3/G_i,w_i)
\rightarrow B_r({\bf H}^3/H', w')$, where $B_r$ denotes an $r$--ball.
(Note that the group $H'$ may depend on $i$ since a geometric limit exists
only after taking a subsequence.)

First suppose that $\lambda'$ can be realized by a pleated surface
homotopic to $\psi'$.
Then by Proposition~\ref{5.1}, for any small $\delta >0$, there exists a
train track $\tau$ with a weight system $\omega$ carrying $\lambda'$ and a
continuous map $f \co  S \rightarrow {\bf H}^3/G'$ homotopic to $\psi'$ which
is adapted to a tied neighbourhood $N_\tau$ of $\tau$ such that the total
curvature and the quadratic variation of angle for  $f(\tau,\omega)$ are
less than $\delta$.

For a Kleinian group $H'$ containing $G'$, by composing the covering $q\co 
{\bf H}^3/G' \rightarrow {\bf H}^3/H'$ to $f$, we get a map with the same
property homotopic to $q \circ \psi'$.
We take $r$ and $\epsilon$ so that for any geometric limit $H'$, the
$r$--ball centred at the base point under $w'$ contains the image of $q
\circ f$ and so that if we pull back $q\circ f$ by a
$((1+\epsilon),r)$--approximate isometry and straighten the images of
branches to geodesic arcs, the image of $(\tau ,\omega)$ has the total
curvature and the quadratic variation of angle less than $2\delta$.
Then for $i>i_0$, there exists a map $f_i \co  S \rightarrow {\bf H}^3/G_i$
homotopic to $\phi_i$ which is adapted to $\tau$ such that
$f_i(\tau,\omega)$ has total curvature and quadratic variation of angle
less than $2\delta$.
Again by Proposition \ref{5.1}, this implies that there is a neighbourhood
$U$ of $\lambda'$ in ${\cal ML}(S)$ such that for any weighted simple
closed curve $\gamma$ in $U$, there exist $\nu_U >0$ depending on $U$,
$\eta_\delta >0$, and $t_\delta<1$ depending on $\delta$ such that $\nu_U
\rightarrow 0$ as $U$ gets smaller and $\eta_\delta \rightarrow 0, t_\delta
\rightarrow 1$ as $\delta \rightarrow 0$, and the following holds.
We can homotope $\gamma$ so that $N_\tau \cap \gamma$ corresponds the
weight system $\omega'$ (which may not satisfy the switch condition since
$\gamma$ may not be homotoped into $N_\tau$) whose value at each branch
differs from that of $\omega$ at most $\nu_U$, and the closed geodesic
$\gamma_i^*$ homotopic to $\psi_i(\gamma)$ has a part with length
$t_{\delta}{\rm length}(\gamma_i^*)$ which lies within distance
$\eta_\delta$ from $f^i(\tau \cap \gamma)$.
The same holds for $f$ and the closed geodesic $\gamma^*$ homotopic to
$\psi(\gamma)$.

It follows that there is a positive real number $\zeta$ depending on
$\epsilon, \delta , U$ which goes to $0$ as $\epsilon \rightarrow 0, \delta
\rightarrow 0$ and $U$ gets smaller remaining to be a neighbourhood of
$\lambda'$, such that if $\gamma, \gamma'$ are weighted simple closed
curves in $U$, then $|{\rm length}(\gamma_i^*) -{\rm
length}({\gamma'}^*)|<\zeta$, where ${\gamma'}^*$ is the closed geodesic in
${\bf H}^3/G'$ homotopic to $\psi(\gamma')$.
Since the set of weighted simple closed curves is dense in ${\cal ML}(S)$
and  any realization of measured lamination can be approximated by
realizations of simple closed curves, this  implies our lemma in the case
when $\lambda'$ is realizable by a pleated surface homotopic to $\psi'$.

Next suppose that $\lambda'$ is not realizable by a pleated surface
homotopic to $\psi'$.
This means that $\lambda'$ is an ending lamination of an end of $({\bf
H}^3/G')_0$.
By a result of Thurston in \cite{Th} (see also Lemma~4.4 in \cite{Oh1}), it
follows that $\lambda'$ is maximal and connected.
In this case the alternative (i) of Proposition \ref{5.1} holds.
Hence for any small $\epsilon >0$, there exists a train track $\tau$
carrying $\lambda'$ with weight $\omega$ and a continuous map $f$ homotopic
to $\psi'$ which is adapted to a tied neighbourhood $N_\tau$ of $\tau$,
such that $\lambda'$ can be homotoped so that the length of $f(\tau,
\omega)$ is less than $\delta$.
Then by the same argument as the last paragraph, there exists $i_0$ such
that if $i>i_0$ there exists a map $f_i$ adapted to $N_\tau$ such that
$f_i(\tau, \omega)$ has length less than $2\delta$.

Since $\{\lambda_j\}$ converges to $\lambda'$ and $\lambda'$ is maximal,
$\lambda_j$ is carried by $\tau$ for sufficiently large $j$ with weight
$\omega_j$ whose values at branches are close to those of $\omega$.
Hence there exists $j_0$ such that $f_i(\tau , \omega_j)$ is less than
$3\delta$ if $j >j_0$.
As the length of realization of $\lambda_j$ by  a pleated surface homotopic
to $\psi_i$ is less than that of $f_i(\tau , \omega_j)$, this implies our
lemma in the case when $\lambda'$ cannot be realized by a pleated surface
homotopic to $\psi'$.
\endproof

The following is a well-known result of Thurston appeared in \cite{Th} and
also a corollary of Lemma \ref{intersection} due to Bonahon.
Nevertheless, as its proof is not so straightforward when sequences of 
closed geodesics intersect  Margulis tubes non-trivially, we shall prove
here that Lemma \ref{intersection} implies this lemma.

\begin{lem}
\label{unique ending}
Let $M$ be a hyperbolic 3--manifold.
Let $e$ be a geometrically infinite tame end of the non-cuspidal part $M_0$.
Let $\lambda, \lambda'$ be measured laminations on a frontier component $T$
of a core, which faces $e$.
Suppose that both $\lambda$ and $\lambda'$ are ending laminations of the
end $e$.
Then the supports of $\lambda$ and $\lambda'$ coincide.
\end{lem} 
 
\noindent
\proof
Let $s_j$ and $s_j'$ be simple closed curves on $T$ such that for some
positive real numbers $x_j$ and $y_j$, we have $x_j s_j \rightarrow
\lambda$, $y_j s_j' \rightarrow \mu$ and such that the closed geodesics
$s_j^*$ homotopic to $s_j$ and ${s_j'}^*$ homotopic to $s_j'$ tend to the
end $e$ as $j \rightarrow \infty$.
If  there exists a constant $\epsilon_0 >0$ such that neither  $s_j^*$ nor
${s_j'}^*$ intersects an $\epsilon_0$--Margulis tube whose axis is not
$s^*_j$ or ${s_j'}^*$ itself, then we can apply Lemma \ref{intersection}
and the proof is completed.

Next suppose that for at least one of $s_j^*$ and ${s_j'}^*$ (say $s_j^*$),
a constant as $\epsilon_0$ above does not exist.
We shall prove that we can replace $s_j$ with another simple closed curve
to which we can apply Lemma \ref{intersection}.
By assumption, there exist closed geodesics $\xi_j$ whose lengths go to $0$
and such that $s_j^*$ intersect the $\epsilon_j$--Margulis tube whose axis
is  $\xi_j$, where $\epsilon_j \rightarrow 0$.
Let $h_j \co  (T,\sigma_j) \rightarrow M$ be a pleated surface homotopic to
the inclusion whose image contains $s_j^*$ as the image of its pleating
locus, where $\sigma_j$ is the hyperbolic structure on $T$ induced by
$h_j$.
Put a base point $y_i$ on $T$ which is mapped into $s_j^*$ but outside  the
$\epsilon_0$--Margulis tubes by $h_j$.
Let $h_\infty \co  ((T',\sigma_\infty),y_\infty) \rightarrow (M',y_\infty)$ be
the geometric limit of $\{h_j\co  (T',y_j) \rightarrow (M,h_j(y_j))\}$  after
taking a subsequence, where $T'$ is an incompressible subsurface in $T$.
We shall first show that $T'$ cannot be the entire of $T$.

Suppose that $T'=T$ on the contrary.
Let $l$ be the geodesic lamination on $(T, \sigma_\infty)$ which is the
geometric limit of the closed geodesic on $(T,\sigma_j)$ corresponding to
$s_j^*$ as $j \rightarrow \infty$.
Since $l$ cannot approach to a cusp (as $T=T'$), it is compact.
Therefore we can take a point in the intersection of $s_j^*$ and the
$\epsilon_j$--Margulis tube which converges to a point $x$ on $h_\infty(l)$
associated with the geometric convergence of $\{h_j\}$ to $h_\infty$ as $j
\rightarrow \infty$.
Then for any small $\epsilon$, there is an essential closed curve passing
$x$ with length less than $\epsilon$ which can be obtained by pushing
forward by an approximate isometry an essential loop intersecting $s_j^*$
of length less than $\epsilon_j$ for sufficiently large $j$.
This is a contradiction.

Thus there is an extra cusp for $h_\infty$.
Let $c$ be a simple closed curve on $T'$ representing an extra cusp.
Let $\overline{\rho}_j \co  B_{r_j}((T, \sigma_j),y_j) \rightarrow
B_{r_j}((T',\sigma_\infty),y_\infty)$ be an approximate isometry associated
with the geometric convergence of $\{h_j\}$ to $h_\infty$.
Let $c_j$ be a simple closed curve on $T$ which is homotopic to
$\overline{\rho}_j^{-1}(c)$.
Let $l'$ be a measured lamination to which $\{r_j c_j\}$ converges for some
positive real numbers $r_j$. 
Let $c_j^+$ be the closed geodesic on $(T,\sigma_j)$ homotopic to $c_j$.
Let $\alpha$ be a measured lamination to whose projective class the
hyperbolic structures $\sigma_j$ converge, after passing through a
subsequence if necessary.
Then as ${\rm length}_{\sigma_j}(c_j^+)$ goes to $0$ as $j \rightarrow
\infty$, we have $i(\alpha, l') =0$ by Lemma~3.4 in \cite{Oh1}.
By the same reason, considering $\{s_j\}$, we have $i(\lambda, \alpha)=0$.
Since $\lambda$ is maximal and connected, these imply that the supports of
$\lambda$ and $l'$ coincide.
In particular, $l'$ is an ending lamination of the end for which $\lambda$
is an ending lamination.

Because the length of  the closed geodesic $c_j^+$  goes to $0$ as $j
\rightarrow \infty$, the closed geodesic homotopic to $h_j(c_j)$, whose
length is at most the length of $c_j^+$, must be the axis of an
$\epsilon_0$--Margulis tube for sufficiently large $j$.
Thus we can replace $s_j$ with $c_j$, and and by the same fashion, we can
replace $s_j'$ with another simple closed curve if necessary.
 We can apply Lemma \ref{intersection} for such simple closed curves.
\endproof

\section{Proof of the main theorem}

We shall complete the proof of Theorem \ref{main} in this section.
Recall that under the assumption for the reductio ad absurdum, we have
$(G,\psi)$ which is the algebraic limit of $\{(\Gamma_i,\phi_i)\}$.

\begin{lem}
\label{ending}
In the situation of Theorem \ref{main}, the non-cuspidal part $({\bf
H}^3/G)_0$ of the hyperbolic 3--manifold ${\bf H}^3/G$ has a geometrically
infinite tame end for which $\psi(\mu)$ represents an ending lamination.
\end{lem}

\noindent
\proof
Suppose first that the end $e^i$ of $({\bf H}^3/\Gamma_i)_0$ corresponding
to the first factor of the end invariant is geometrically finite.
Then by Lemma~\ref{sequence}, there exists a sequence of weighted simple
closed curves $r_i\gamma_i$ on $S$ converging to $\mu$ such that for the
closed geodesic $\gamma_i^*$ in ${\bf H}^3/\Gamma_i$ homotopic to
$\phi_i(\gamma_i)$, we have $r_i {\rm length}(\gamma_i^*) \rightarrow 0$.
By the continuity of length function $L$ on $AH_p(S) \times {\cal ML}(S)$
(Lemma~\ref{continuity}), we have $L((G,\psi),\mu)$ is $0$, which means
that $\mu$ cannot be realized by a pleated surface homotopic to $\psi$.
As we assumed that $\mu$ is maximal and connected, there must be a
geometrically infinite tame end of $({\bf H}^3/G)_0$ with ending lamination
represented by $\psi(\mu)$.
This last fact, originally due to Thurston, can be proved using Bonahon's
result:
by  Proposition \ref{5.1}, if $L((G,\psi),\mu)=0$ and $\mu$ is maximal and
connected, then for any sequence of simple closed curves $\delta_j$ on $S$
whose projective classes converge to that of $\mu$, the closed geodesics
$\delta_j^*$ homotopic to $\psi(\delta_j)$ tend to an end of $({\bf
H}^3/G)_0$.
This means that $\psi(\mu)$ is an ending lamination for a geometrically
tame end of $({\bf H}^3/G)_0$.

Next suppose that the end $e^i$ is geometrically infinite.
Then $\chi_i$ is represented by a measured lamination $\mu_i$ which
represents an ending lamination of $e^i$, hence $L((\Gamma_i,\phi_i),
\mu_i)=0$.
We can assume that $\mu_i$ lies on the unit ball of ${\cal ML}(S)$ with
respect to the metric induced from some fixed hyperbolic structure on $S$.
Then $\mu_i$ converges to a scalar multiple of $\mu$ since we assumed that
$\chi_i=[\mu_i]$ converges to $[\mu]$.
By the continuity of $L$, this implies that $L((G,\psi), \mu)=0$ and that
$\psi(\mu)$ represents  an ending lamination  for $({\bf H}^3/G)_0$.
\endproof

We shall denote the end in Lemma~\ref{ending}, for which $\psi(\mu)$
represents an ending lamination, by $e_\mu$.

Recall that $q \co {\bf H}^3/G \rightarrow {\bf H}^3/\Gamma_\infty$ is
a covering associated with the inclusion.  Now by Proposition
\ref{covering}, the end $e_\mu$ has a neighbourhood $E_\mu$ such that
$q|E_\mu$ is a proper embedding.  Since $\mu$ is maximal and
connected, the end $e_\mu$ has a neighbourhood homeomorphic to $S_0
\times {\bf R}$, where $S_0$ is the non-cuspidal part of $S$ Hence by
refining $E_\mu$, we can assume that $E_\mu$ is also homeomorphic to
$S_0 \times {\bf R}$.

\begin{lem}
\label{far surface}
We can take an embedding $f' \co  S_0 \rightarrow E_\mu$ homotopic to
$\psi|S_0$ whose image is contained in the convex core such that for any
pleated surface $g\co  S \rightarrow {\bf H}^3/\Gamma_\infty$ homotopic to $q
\circ \psi$ with non-empty intersection with $qf'(S_0)$, we have $g(S) \cap
({\bf H}^3/\Gamma_\infty)_0  \subset q(E_\mu)$.
\end{lem}

\noindent
\proof
Fix a constant $\epsilon_0>0$ less than the Margulis constant.
There exists a constant $K$ such that for any hyperbolic metric on $S$, the
diameter of $S$ modulo the $\epsilon_0$--thin part is bounded above by $K$.
(This can be easily seen by considering the moduli of $S$.)

Note that since the end $e_\mu$ is geometrically infinite, it has a
neighbourhood contained in the convex core.
Take $t \in {\bf R}$ large enough so that $S_0 \times \{t\}
\subset E_\mu$ is contained in the convex core 
and the distance from  $S_0 \times \{t\}$ to the
frontier of $E_\mu$ in $({\bf H}^3/G)_0$  modulo the $\epsilon_0$--thin part
is greater than $2K$.
Choose $f'$ homotopic $\psi|S_0$ so that its image is $S_0 \times \{t\}$. 
Then the distance between $q f'(S_0)$ and the frontier of $q(E_\mu)$ modulo
the $\epsilon_0$--thin part is also greater than $2K$.
Suppose that a pleated surface $g \co  S \rightarrow {\bf H}^3/\Gamma_\infty$
touches $qf'(S_0)$.
Then $g(S)$ cannot meet the frontier of $q(E_\mu)$ since the
$\epsilon_0$--thin part of $S'$ with respect to the hyperbolic structure
induced by $g$ is mapped into the $\epsilon_0$--thin  part of ${\bf
H}^3/\Gamma_\infty$, and any path on $g(S)$ has length less than  $K$
modulo the $\epsilon_0$--thin part of ${\bf H}^3/\Gamma_\infty$.
Also it is impossible for $g(S)$ to go into the cuspidal part of ${\bf
H}^3/\Gamma_\infty$ and come back to the non-cuspidal part since the
intersection of $g(S)$ with the cuspidal part of ${\bf H}^3/\Gamma_\infty$
is contained in a neighbourhood of cusps of $g(S)$.
This means that without meeting the frontier of $q(E_\mu)$, the pleated
surface $g(S)$ cannot  go outside $q(E_\mu)$ in $({\bf
H}^3/\Gamma_\infty)_0$.
Thus the intersection of such a pleated surface with $({\bf
H}^3/\Gamma_\infty)_0$ must be contained in $q(E_\mu)$.
\endproof

We denote $q \circ f'$ by  $f_\infty \co  S_0 \rightarrow ({\bf
H}^3/\Gamma_\infty)_0$.
Pulling back this embedding $f_\infty$ by an approximate isometry $\rho_i$
for sufficiently large $i$, we get an embedding $f_i \co  S_0 \rightarrow
({\bf H}^3/\Gamma_i)_0$.
Since $f_\infty$ comes from the surface homotopic to $\psi$ in the
algebraic limit, for sufficiently large $i$, the surface $f_i$ is homotopic
to $\phi_i$.

Consider the case when $({\bf H}^3/\Gamma_i)_0$ has a geometrically finite
end; that is $\Gamma_i$ is either quasi-Fuchsian or a totally degenerate
b--group.
As in the previous section, let $\Sigma_i$ be a  boundary components of the
convex core of ${\bf H}^3/\Gamma_i$, and let $h_i \co  S \rightarrow \Sigma_i$
 be a homeomorphism homotopic to $\phi_i $.
The homeomorphisms $h_i$ can also be regarded as pleated surfaces in ${\bf
H}^3/\Gamma_i$.
Let $\mu_i$  be the bending locus of $h_i$, to which we give transverse
measures with full support so that $\mu_i$  should converge to measured
laminations as $i \rightarrow \infty$  after taking subsequences.
(Since the unit sphere of the measured lamination space is compact, this is
always possible.
Also if $h_i$ happens to be totally geodesic, we can set $\mu_i$ to be  any
measured lamination on $S$.)

\begin{lem}
\label{pleats limit}
Suppose that $\Gamma_i$ is either quasi-Fuchsian or a totally degenerate
{\rm b}--group as above.
The sequence of the measured laminations $\{\mu_i\}$  converges to a
measured lamination with the same support as $\mu$ after taking
subsequences.
\end{lem}

\noindent
\proof
Let $\mu'$ be a limit of $\{\mu_i\}$ after taking a subsequence.
If $i(\mu,\mu')=0$, we have nothing to prove any more because $\mu$ is
maximal and connected.
Now assume that $i(\mu',\mu) \neq 0$.
Then, by the fact that the marked hyperbolic structure on $\Sigma_i$
converges to $[\mu]$ as $i \rightarrow \infty$ and Lemma~3.4 in \cite{Oh1},
we have ${\rm length}_{\Sigma_i}(\mu_i) \rightarrow \infty$.
On the other hand, by the continuity of the length function $L$ on $AH_p(S)
\times {\cal ML}(S)$ (Lemma~\ref{continuity}), we have $${\rm length}_{{\bf
H}^3/\Gamma_i}(\phi_i(\mu_i))=L((\Gamma_i,\phi_i),\mu_i) \rightarrow
L((G,\psi),\mu')={\rm length}_{{\bf H}^3/G}(\psi(\mu')) < \infty$$ where
${\rm length}_{{\bf H}^3/\Gamma_i} (\phi_i(\mu_i))$ denotes the length of
the image of $\mu_i$ realized by pleated surface homotopic to $\phi_i$ etc.
Since $\mu_i$ is mapped by $h_i$ into the bending locus of $\Sigma_i$, it
is realized by $h_i$, hence ${\rm length}_{\Sigma_i}(\mu_i)=
{\rm length}_{{\bf H}^3/\Gamma_i}(\phi_i(\mu_i))$.
This is a contradiction.
\endproof

Now we assume that $\Gamma_i$ is quasi-Fuchsian.
Then there are two boundary components $\Sigma_i, \Sigma_i'$ of the convex
core of ${\bf H}^3/\Gamma_i$, and homeomorphisms $h_i \co  S \rightarrow
\Sigma_i \subset {\bf H}^3/\Gamma_i$ and $h_i' \co  S\rightarrow \Sigma_i'
\subset {\bf H}^3/\Gamma_i$ homotopic to $\phi_i$ which are regarded as 
pleated surfaces.
We have two measured laminations of unit length $\mu_i$ and $\mu_i'$  whose
supports are the bending loci of $h_i$ and $h_i'$.
By Lemma~\ref{pleats limit}, the sequence of the measured laminations
$\{\mu_i\}$ converges to a measured lamination $\mu'$ and $\{\mu_i'\}$
converges to a measured lamination $\mu''$ such that
$|\mu'|=|\mu''|=|\mu|$.
As the space of transverse measures on a geodesic lamination is connected
(or more strongly, convex with respect to the natural PL structure), we can
join $\mu'$ and $\mu''$ by an arc $\alpha \co  I=[0,1] \rightarrow {\cal
ML}(S)$ such that $|\alpha(t)|=|\mu|$.
Join $\mu_i$ and $\mu_i'$ by an arc $\alpha_i \co  I \rightarrow {\cal ML}(S)$
which converges to the arc $\alpha$ joining $\mu'$ and $\mu''$.

Next suppose that $\Gamma_i$ is a totally degenerate b--group.
We can assume without loss of generality that the first factor $\chi_i$ of
the end invariant represents an ending lamination and the second
$\upsilon_i$ a conformal structure. 
Then  we have a pleated surface $h_i \co  S \rightarrow \Sigma_i$ homotopic to
$\phi_i$ whose image is  the boundary  of the convex core.
Let $\mu_i$ be a measured lamination of the unit length whose support is
equal to that of the bending locus as before.
Again by Lemma \ref{pleats limit}, we see that $\{\mu_i\}$ converges to a
measured lamination $\mu'$ with the same support as $\mu$. 
Let $\mu'_i$ be a measured lamination of the unit length representing the
class $\chi_i$.
By the assumption of Theorem \ref{main}, the sequence $\{\mu_i'\}$
converges to a measured lamination $\mu''$ with the same support as $\mu$.
As in the case of quasi-Fuchsian group, we join $\mu'$ and $\mu''$ by an
arc $\alpha$, and then  join $\mu_i$ and $\mu_i'$ by an arc $\alpha_i$
which does not pass an ending lamination for ${\bf H}^3/\Gamma_i$ at the
interior so that it will converge to $\alpha$ uniformly.

In the case when $\Gamma_i$ is a totally doubly degenerate group, both
$\chi_i$ and $\upsilon_i$ are represented by  ending laminations.
Let $\mu_i$ representing $\chi_i$ and $\mu_i'$ representing $\upsilon_i$ be
measured laminations of the unit length.
Then by assumption, $\{\mu_i\}$ and $\{\mu_i'\}$ converge to measured
laminations $\mu'$ and $\mu''$ with the same support as $\mu$.
As before we join $\mu'$ and $\mu''$ by an arc $\alpha$, and $\mu_i,
\mu_i'$ by $\alpha_i$ which does not pass an ending lamination of ${\bf
H}^3/\Gamma_i$ at the interior so that $\{\alpha_i\}$ converges to
$\alpha$.

Next we shall consider constructing for each $i$ a homotopy consisting of
pleated surfaces and negatively curved surfaces in ${\bf H}^3/\Gamma_i$ as
in Thurston~\cite{Th}.
What we shall have is a homotopy between $h_i$ and $h_i'$ in the case when
$\Gamma_i$ is quasi-Fuchsian; a half-open homotopy $\hat{H}_i \co  S \times
[0,1) \rightarrow {\bf H}^3/\Gamma_i$ such that $\hat{H}_i(S \times \{t\})$
tends to the unique geometrically infinite end as $r \rightarrow 1$ in the
case when $\Gamma_i$ is a totally degenerate b--group; and an open homotopy
$\hat{H}_i \co  S \times (0,1) \rightarrow {\bf H}^3/\Gamma_i$ such that
$\hat{H}_i(S\times \{t\})$ tends to one end as $ t\rightarrow 0$ and to the
other as $t \rightarrow 1$ in the case when $\Gamma_i$ is a totally doubly
degenerate group.
To construct such a homotopy, we need the notion of rational depth for
measured laminations due to Thurston.
An alternative approach to construct such a homotopy using  singular
hyperbolic triangulated surfaces can be found in Canary~\cite{Ca}.

A train track $\tau$ is called birecurrent when the following two
conditions are satisfied.
(This definition is due to Penner--Harer~\cite{PH}.)
(1) The $\tau$ supports a  weight system which is positive on each branch
$b$ of $\tau$.
(2) For each branch $b$ of $\tau$, there exists a multiple curve $\sigma$
(ie a disjoint union of non-homotopic essential simple closed curves)
transverse to $\tau$ which intersects $b$ such that $S-\tau -\sigma$ has no
bigon component.\\
A birecurrent train track which is not a proper sub-train track of another
birecurrent train track is said to be complete.

Any measured lamination is carried by some complete train track.
(Refer to Corollary 1.7.6 in \cite{PH}.)
The weight systems on a complete train track gives rise to a coordinate
system in the measured lamination space.
(See Lemma 3.1.2 in \cite{PH}.)
The rational depth of a measured lamination is defined to be the dimension
of the rational vector space of linear functions with rational coefficients
vanishing on the measured lamination with respect to a coordinate system
associated with a complete train track carrying the measured lamination.
This definition is independent of the choice of a coordinate system since
functions corresponding to coordinate changes are linear functions with
rational coefficients.
The set of measured laminations with rational depth $n$ has codimension $n$
locally at regular points.
In particular a generic arc in the measured lamination space does not pass
a measured lamination with rational depth more than $1$.

Now perturb $\alpha$ and $\alpha_i$ to a piece-wise linear path with
respect to the PL structure of ${\cal ML}(S)$ determined by complete train
tracks fixing the endpoints so that for each $t \in I$, the measured
lamination $\alpha_i(t)$ is not an ending lamination and has rational depth
$0$ or $1$,  and that for each $i$  there exist only countably many values
$t$ for which $\alpha_i(t)$ has rational depth 1.

The following lemma was first proved in section 9 in Thurston~\cite{Th}.
A fairly detailed proof can be found there.

\begin{lem}
\label{depth 0}
If a measured lamination has rational depth $0$, then each component of its
complement is either an ideal triangle  or a once-punctured monogon except
when $S$ is a once-punctured torus.
In the case when $S$ is a once-punctured torus, each component  of the
complement is an ideal once-punctured bigon.
A pleated surface $f \co  S \rightarrow M$ realizing a measured lamination
$\zeta$ of rational depth $0$ is unique among the maps in the homotopy
class, and every sequence of homotopic pleated surfaces realizing measured
laminations converging to $\zeta$ converges to the pleated surface
realizing $\zeta$.
\end{lem}

\noindent
\proof
First we shall show that each complementary region of a measured lamination
of rational depth $0$ is either an ideal triangle or an ideal 
once-punctured monogon unless $S$ is a once-punctured torus.

Suppose that $S$ is not a once-punctured torus and that a measured
lamination $\zeta$ has a complementary region which is neither  an ideal
triangle nor an ideal once-punctured monogon.
Then, we can construct a birecurrent train track $\tau$ carrying $\zeta$
whose complement has a component which is neither a triangle nor a
once-punctured monogon.
(Refer to section 1.7 in \cite{PH}.)

A birecurrent train track is maximal if and only if every component of its
complement is either a triangle or a once-punctured monogon, and
non-maximal birecurrent train track is a sub-train track of a complete
train track.
(Theorem 1.3.6 in \cite{PH}.)
Hence there exists a complete train track $\tau'$ containing $\tau$ as a
proper sub-train track.
Since there is a branch of $\tau'$ through which $\zeta$ does not pass
after homotoping $\zeta$ so that it is carried by $\tau'$, it follows that
with respect to the coordinate system corresponding to $\tau'$, the
measured lamination $\zeta$ has rational depth at least $1$.

In the case when $S$ is a once-punctured torus, again Theorem 1.3.6 in
\cite{PH} says that a birecurrent train track is maximal if and only if its
(unique) complementary region is a once-punctured bigon.
Thus the same argument as above also implies our claim in the case of
once-punctured torus.

Next we shall show the uniqueness of realization of a measured lamination
of rational depth $0$.
Let $f, g$ be two pleated surfaces realizing a measured lamination $\zeta$
of depth $0$.
The pleated surfaces $f,g$  induce hyperbolic metrics $m_1, m_2$
respectively  on $S$.
(These may differ as we do not know if $f$ and $g$ coincide.)
The measured lamination $\zeta$ is homotopic to  measured geodesic
laminations $\zeta_1$ with respect to $m_1$ and $\zeta_2$ with respect to
$m_2$.
Consider the universal covers $p_1 \co  {\bf H}^2 \rightarrow (S,m_1)$ and
$p_2 \co  {\bf H}^2 \rightarrow (S,m_2)$.
Let $\tilde{\zeta}_1$ be $p_1^{-1}(\zeta_1)$ and let $\tilde{\zeta}_2$ be 
$p_2^{-1}(\zeta_2)$.

The pleated surfaces $f,g$ are lifted to maps $\tilde{f}, \tilde{g} \co  {\bf
H}^2 \rightarrow {\bf H}^3$.
Since $\zeta$ has compact support, there is a homeomorphism from $S$ to $S$
homotopic to the identity which takes $\zeta_1$ to $\zeta_2$ and is equal
to the identity near cusps. 
Also for a homotopy between $f$ and $g$, the distance moved by the homotopy
on the compact set $\zeta$ has an upper bound. 
These imply that for each leaf $l$ of $\tilde{\zeta}$ the images of the
corresponding leaves $l_1$ of $\tilde{\zeta}_1$ by $\tilde{f}$ and $l_2$ of
$\tilde{\zeta}_2$ by $\tilde{g}$ are within a bounded distance. 
Since both $\tilde{f}(l_1)$ and $\tilde{g}(l_2)$ are geodesics in ${\bf
H}^3$ and two geodesics lying in bounded distance coincides in ${\bf H}^3$,
these two images must coincide.
Hence we have a map $q \co  {\bf H}^2 \rightarrow {\bf H}^2$ equivariant with
respect to the action of $\pi_1(S)$ with the property $\tilde{f}|
\tilde{\zeta}_1 = \tilde{g}\circ q| \tilde{\zeta}_1$ which maps
$\tilde{\zeta}_1$ to $\tilde{\zeta}_2$ isometrically.

It remains to prove that $q$ extends to an equivariant  isometry
$\overline{q}$ of ${\bf H}^2$ with  the property $\tilde{f} = \tilde{g}
\circ \overline{q}$.
Since $\zeta$ has rational depth $0$, each of its complementary regions is
either an ideal triangle or an ideal once-punctured monogon unless $S$ is a
once-punctured torus.
An ideal triangle on $S$ is lifted to that on ${\bf H}^2$.  Since the three
sides of the triangle are mapped to geodesics by $\tilde{f}$ or
$\tilde{g}$, the triangle must be mapped totally geodesically.
Considering that there is only one isometry type of ideal triangles, we can
see that this implies $q$ can be extended to ideal triangle complementary
components  without problem.

For complementary regions which are ideal once-punctured monogon, or ideal
once-punctured bigon in the case when $S$ is a once-punctured torus, we
have to use the fact that pleated surfaces are totally geodesics near
cusps.
(This is proved in Proposition 9.5.5 in Thurston \cite{Th}.)
Once we know this, we can subdivide such regions into ideal triangles by
adding geodesics tending to cusps on $S$, which are mapped to geodesics by
$f$ or $g$.
Since each cusp of $S$ is mapped to the same cusp of $M$ by $f$ and $g$, we
can arrange them so that the lifts of these added geodesics should be compatible
with $q$.
Hence by extending the map finally to ideal triangles, we get a map
$\overline{q}$ as we wanted.

Finally let us prove the last sentence of our lemma.
Let $\xi_j$ be measured laminations converging to $\zeta$, and $f_j$ a
pleated surface realizing $\xi_j$.
Since $\zeta$ can be realized by a pleated surface, the alternative (2) of
Proposition~\ref{5.1} should be valid for $\zeta$.
We shall show that if there is no compact set in $M$ which intersects all
the images of $f_j$, then we can  see that the alternative (2) of 
Proposition \ref{5.1} fails to hold for $\zeta$.

Suppose that the alternative (2) of Proposition~\ref{5.1} holds for $\zeta$.
Then for  any $\delta >0$ and $t <1$, there exist a map $f_\delta \co  S
\rightarrow M$ homotopic to $f$ such that for any simple closed curve
$\gamma$ whose projective class is close to that of $\zeta$,  the closed
geodesic $\gamma^*$ homotopic to $f_\delta(\gamma)$ has a part of length at
least $t{\rm length}(\gamma^*)$ which is contained in the
$\delta$--neighbourhood of $f_\delta(\zeta)$.
Note that as $\delta \rightarrow 0$, this map $f_\delta$ converges to a
pleated surface realizing $\zeta$, which must be equal to $f$.
(Refer to \cite{Oh4} for a further explanation.)
On the other hand, since $\xi_j$ is also realized by a pleated surface homotopic to $f$, the alternative (2) holds also for $\xi_j$.
Then we have a surface $f_j^\delta$ with the same property as $f_\delta$
above replacing $\zeta$ with $\zeta_j$.
Since we assumed that $f_j$ tends to an end of $M$, we can have  surfaces
$f_j^{\delta_j}$ going to an end and a simple closed curve $\gamma_j$ whose
projective class is close to that of $\zeta_j$ such that a large part of
the closed geodesic $\gamma_j^*$ is contained in the
$\delta_j$--neighbourhood of $f_j^{\delta_j}(S)$.
This is a contradiction because $\gamma_j^*$ must also have a large part
contained in the $\delta$--neighbourhood of $f_\delta(S)$ which remains in a
neighbourhood of $f(S)$.

Thus the surfaces $f_j(S)$ remain to intersect a compact set, hence
converge to a pleated surface $g$ homotopic to $f$ uniformly on any compact
set of $S$.
(Theorem~5.2.18 in Canary--Epstein--Green~\cite{CEG}.)
The pleated surface $g$ realizes a geodesic lamination $\zeta_\infty$ which
is a geometric limit of $\{\zeta_j\}$ regarded as geodesic laminations
forgetting the transverse measures.
It is known that $\zeta_\infty$ contains the support of $\zeta$.
(See for example Lemma 5.3.2. in \cite{CEG}.)
Thus $g$ also realizes $\zeta$, and by the uniqueness of such pleated
surfaces proved above, we see that $f$=$g$, which means that $\{f_j\}$
converges to $f$ uniformly on any compact set of $S$.
\endproof

Now let $H_i \co  S \times I, S \times [0,1), S \times (0,1) \rightarrow {\bf
H}^3/\Gamma_i$ (depending on the type of $\Gamma_i$; a quasi-Fuchsian group
or  a totally degenerate b--group or  a totally doubly degenerate) be a map
such that for each $t \in I, [0,1), (0,1) $, the map $H_i(\ ,t) \co  S
\rightarrow {\bf H}^3/\Gamma_i$ is a pleated surface realizing
$\alpha_i(t)$.
Then $H_i$ is continuous with respect to $t$ by Lemma~\ref{depth 0} except
at values $t$ where $\alpha(t)$ has rational depth $1$, which are countably
many.
Since we made $\alpha_j$ piece-wise linear we can see that the right and
left limits exists even at $t$ where $\alpha_i(t)$ has depth $1$.
(This can be seen by considering a complete train track giving a coordinate
near $\alpha_i(t)$.)
As was shown in Thurston~\cite{Th}, (see section~4 in \cite{Oh4} for an
explanation), at such a point of discontinuity $t$, the left limit and the
right limit differ only in a complementary region $R$ of $\alpha_i(t)$
which is either an ideal quadrilateral or an ideal once-punctured bigon
except when $S$ is a once-punctured torus.
Then we can modify $H_i$ to a continuous homotopy $\hat{H}_i$ by
interpolating negatively curved surfaces realizing $\alpha_i(t_0)$ between
$\lim_{t\rightarrow t_0-0}H_i(\ ,t)$ and $\lim_{t\rightarrow t_0+0}H_i(\
,t)$ at each $t_0$ where $\alpha_i(t_0)$ has rational depth $1$ as in
Thurston~\cite{Th}.
These negatively curved surface coincide with the left and the right limit
outside $R$ where the right and left limit differ.

We need to prove that a family of surfaces thus obtained is continuous with
respect to the parameter.
The only case that we have to take care of is when the values $t_k$ for
which $\alpha_i(t_k)$ has depth $1$ accumulates to a point $t_0 \in I$.
The negatively curved surfaces interpolated at $t_k$ have the same image as
a pleated surface realizing $\alpha(t_k)$ outside a complementary region
$R_k$.
The image of $R_k$ by the left limit pleated surface and the right limit
surface  bound an ideal tetrahedron if $R_k$ is an ideal quadrilateral or a
solid torus with cusps if $R_k$ is an ideal once-punctured bigon in the
case when $S$ is not a once-punctured torus.
The form of $R_k$ gets thinner and thinner as $k \rightarrow \infty$ since
$t_k$ accumulates.
(This can again be seen by considering a coordinate chart given by a
complete train track.)
This implies that the trajectories  of the homotopy between the left limit
and the right limit, which are contained in the ideal tetrahedron or the
solid torus have length going to $0$ as $k \rightarrow \infty$.
Even in the case when $S$ is a once-punctured torus, a similar argument can
work although we need to take more cases into account.
Thus we can see that $\hat{H}(\ , t)$ is continuous with respect $t$ even
at the point $t_0$ to which depth-1 points $t_k$ accumulate.

\begin{lem}
\label{essential}
For each $i$, there is a pleated surface $k_i\co  S \rightarrow {\bf H}^3/G$
homotopic to $\phi_i$ touching $f_i(S_0)$ which realizes a measured
lamination $\overline{\mu}_i$ such that $\{\overline{\mu}_i\}$ converges to
a measured lamination $\overline{\mu}$ with the same support as $\mu$ after
taking a subsequence.
Moreover we can choose $k_i$ so that $k_i^{-1}(f_i(S_0))$ contains an
essential component relative to cusps.
\end{lem}

\begin{rem}
Although the last sentence of this lemma is not necessary for our purpose
now,  it will be used for our forthcoming work in \cite{Ohf}.
Also Canary's result on filling a convex core by pleated surfaces in
\cite{Ca} will  suffice to prove only the former part of our lemma.
\end{rem}
 
\noindent
\proof
Recall that $f'(S_0)$ is contained in the convex core of ${\bf H}^3/G$.
Then we can assume that $f_i(S_0)$ is also contained in the convex core of
${\bf H}^3/\Gamma_i$.
It follows that $f_i(S_0)$ is contained in the image of $\hat{H}_i$.

By perturbing $f_i(S_0)$ if necessary, we can assume that $\hat{H}_i$ is
transverse to $f_i(S_0)$ and that $\hat{H}_i^{-1}(f_i(S_0))$ is an embedded
surface in $S \times (0,1)$.
Let $F$ be a component of $\hat{H}_i^{-1}(f_i(S_0))$ which separates $S
\times \{0\}$ from $S\times\{1\}$.
It is easy to see such a component exists because in the case when
$\Gamma_i$ is quasi-Fuchsian, $\Sigma_i$ and $\Sigma_i'$ lie in different
components  of ${\bf H}^3/\Gamma_i - f_i(S_0)$, in the case when $\Gamma_i$
is a b--group, $f_i(S_0)$ separates a geometrically infinite end from
$\Sigma_i$, and in the case when $\Gamma_i$ is doubly degenerate,
$f_i(S_0)$ separates two ends.
Then $\pi_1(F)$ is mapped  onto $\pi_1(S)$ by the homomorphism induced by
inclusion, hence $(\hat{H}_i|F)_\# \co  \pi_1(F) \rightarrow \pi_1(f_i(S_0))$
is surjective.

We can assume that for each $t\in I$, the intersection $(S \times
\{t\})\cap F$ is at most one dimensional by perturbing $f_i(S_0)$ again if
necessary.
Then there exists $t_0 \in I$ such that $(S \times \{t_0\}) \cap F$
contains a simple closed curve $K$ which represents a non-trivial element
of $\pi_1(S)$ relatively to the punctures of $S$.
If $\hat{H}_i(\ ,t_0)$ is a pleated surface, we simply let  $k_i$ be
$\hat{H}_i(\ , t_0)$.
In this case, $k_i^{-1}(f_i(S_0))$ contains $K$, which is  essential
relatively to the cusps.
The pleated surface $k_i$ realizes a measured lamination $\alpha_i^{t_0}$
in the image of $\alpha_i$, which we let be $\overline{\mu}_i$.
The measured lamination $\overline{\mu}_i=\alpha_i^{t_0}$ converges after
taking a subsequence to a measured lamination in $\alpha(I)$ hence with the
same support as $\mu$.

Suppose that $\hat{H}_i(\ , t_0)$ is an interpolated negatively curved surface.
Let $\alpha^{t_0}_i$ be the measured lamination of rational depth $1$
realized by $\hat{H}_i(\ ,t_0)$.
We can assume that $\hat{H}_i(\ ,t_0)(\alpha^{t_0}_i)$ is transverse to
$f_i(S_0)$ again by a  perturbation of $f_i(S_0)$ without changing the
homotopy class of $K$.
Let $J=[t_0,t_1]\subset I$ be an interval such that $\hat{H}_i(\
,[t_0,t_1))$ are interpolated negatively curved surfaces and $\hat{H}_i(\
,t_1)$ is a pleated surface realizing $\alpha^{t_0}_i$.
Let $C$ be a component of the complement of $\alpha^{t_0}_i$ which is not
an ideal triangle.
Since $\alpha^{t_0}_i$ has rational depth $1$, such a component is unique
and every simple closed curve in $C$ is either represents a cusp or
homotopic to ${\rm Fr}C$.

On the other hand, by the construction of interpolated surfaces,
$\hat{H}_i|{\rm Fr C} \times J$ is constant with respect to $ t \in J$.
If $C \times \{t_0\}$ does not intersect $K$, the pleated surface
$\hat{H}_i(\ ,t_1) \cap F$ contains a simple closed curve homotopic to $K$,
and we can let $\hat{H}_i(\ ,t_1)$ be $k_i$.
Suppose that $C \times \{t_0\}$ intersects $K$.

\begin{figure}[htb]
\centerline{\relabelbox\small
\epsfxsize 2.5truein\epsfbox{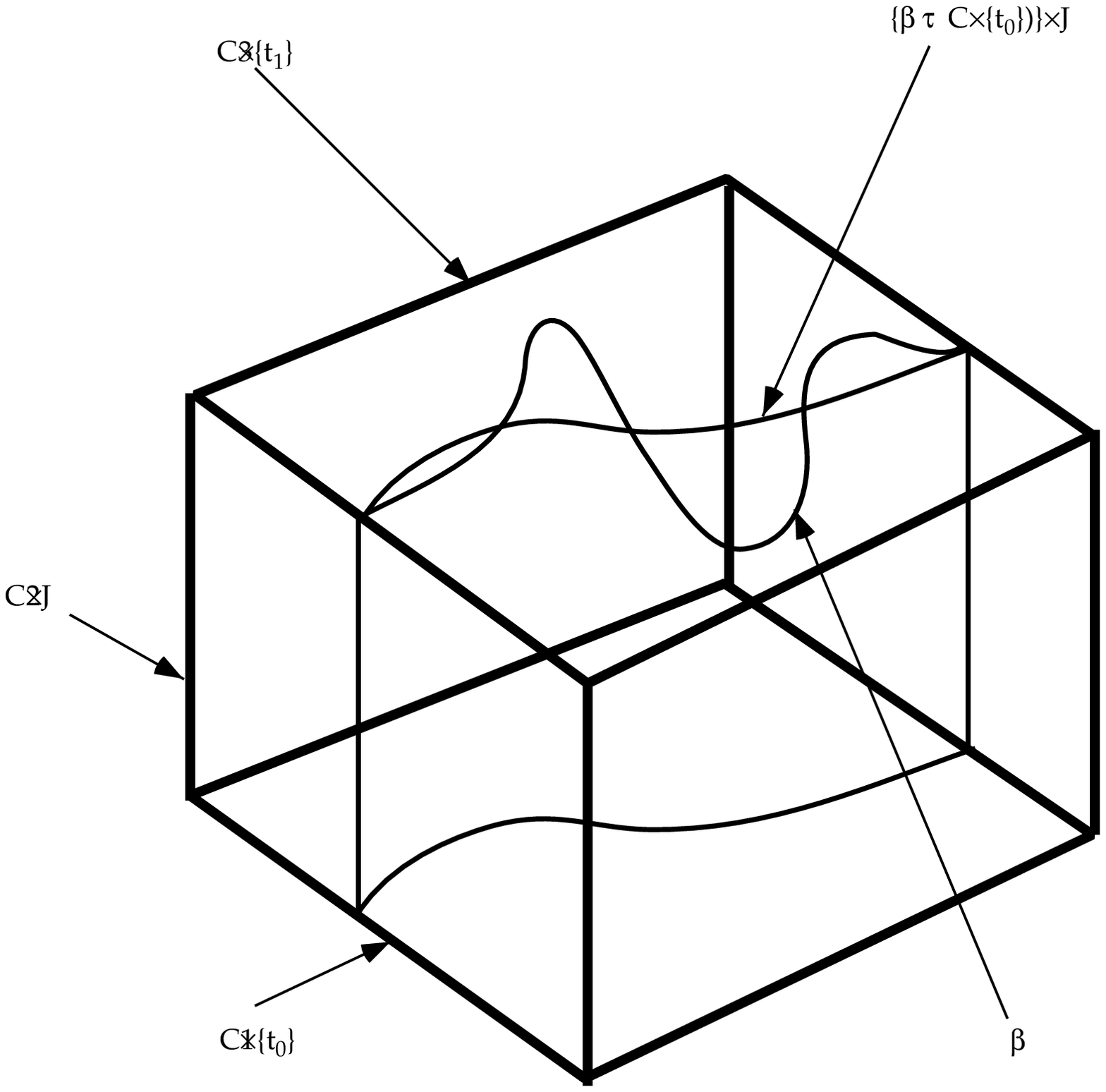}
\adjustrelabel <-15pt, 2pt> {C3}{$C\times\{t_1\}$}
\adjustrelabel <-15pt, 0pt> {C2}{$C\times J$}
\adjustrelabel <-15pt, -2pt> {C1}{$C\times\{t_0\}$}
\adjustrelabel <-2pt, -2pt> {b}{$\beta$}
\adjustrelabel <-30pt, -1pt> {bt}{$\big(\beta\cap(C\times\{t_0\})\big)\times J$}
\endrelabelbox}
\end{figure}

First consider the case when $S$ is not a once-punctured torus.
Then $C$ is either simply connected or an ideal once-punctured monogon.
Consider a component $\beta$ of $(C \times J) \cap F$ intersecting  $K$.
Since $K$ does not represent a cusp, each component of $\beta \cap (C
\times \{t_0\}) \cap K$ must be an open arc.
Since $\hat{H}_i|{\rm Fr}C \times J$ is constant with respect to $t \in J$,
the component $\beta$ must be isotopic to $\{\beta \cap (C \times\{t_0\})
\}\times J$ fixing $\overline{\beta} \cap ({\rm Fr} C \times J)$.
This implies that there exists a component $K'$ of $S \times \{t_1\} \cap
F$ such that $\hat{H}_i(K',t_1)$ is homotopic in ${\bf H}^3/\Gamma_i$ to
$\hat{H}_i(K,t_0)$ on $f_i(S_0)$.
Hence by letting $\hat{H}_i(\ ,t_1)$ be $k_i$, we get a surface as we wanted.

Next suppose that $S$ is a once-punctured torus.
The only case to which the argument above cannot be applied is one when $C$
is a once-punctured open annulus and $K$ is contained in $C \times\{t_0\}$.
By isotoping $f_i(S_0)$ if necessary we can assume that  all the components
of $(C \times J) 
\cap F$ are annuli.
Still there is a possibility that the component of $(C \times J) \cap F$
containing $K$ is an annulus which is parallel into $ C \times \{t_0\}$,
and our argument above would break down.
If there is another essential  (ie incompressible and not
boundary-parallel, where we regard $C \times \partial J$ as the boundary,)
component of $(C \times J) \cap F$, then we can retake $K$ so that $K$ lies
on its boundary and our argument above can be applied.
Suppose that all the components are inessential.
Then consider another interval $J'=[t_2,t_0] \subset I$, such that
$\hat{H}_i(\ ,t)$ is an interpolated surface if $t \in (t_2,t_0]$ and
$\hat{H}_i(\ ,t_2)$ is a pleated surface realizing $\alpha_i^{t_0}$.
Again we can assume that all the components of $(C \times J') \cap F$ are
annuli.
Then some component of $(C \times J') \cap F$ is essential because
otherwise $F$ cannot be a surface separating $S \times \{0\}$ from $S
\times\{1\}$. Hence by the argument as before, retaking $K$, we can assume
that the component of $(S \times J') \cap F$ containing $K$ intersects $S
\times \{t_2\}$ by a simple closed curves homotopic to $K$.

Thus in either case, we can get a pleated surface $k_i$ realizing
$\alpha_i^{t_0}$, which is either $\hat{H}_i(\ , t_1)$ or $\hat{H}_i(\
,t_2)$, and which intersects $f_i(S_0)$ so that the inverse image of
$f_i(S_0)$ has a non-contractible component that is not homotopic to a
cusp.
\endproof

\noindent
\proof[Proof of Theorem \ref{main}]
Consider a geometric limit $k_\infty \co  S' \rightarrow {\bf
H}^3/\Gamma_\infty$ of the sequence of pleated surfaces $k_i \co  S
\rightarrow {\bf H}^3/\Gamma_i$ constructed above.
(Here $S'$ is an open incompressible surface on $S$.)
By construction, $k_\infty(S)$ intersects $f_\infty(S_0)$.
Suppose that $S'$ is not equal to $S$.
Then there is a frontier component $c$ of $S'$ on $S$ which does not
represent a cusp of $S$.
Note that we can apply the same argument as Lemma \ref{far surface}, and
prove that $k_\infty(S)$ does not meet the frontier of $q(E_\mu)$.
Now since $k_\infty(c)$ is homotopic to a cusp component  of ${\bf
H}^3/\Gamma_\infty$ which can be reached from $q(E_\mu)$, it is homotopic
to the image of a cusp of $S$ by $f_\infty$.
By pulling back a homotopy by an approximate isometry, this implies that
$k_i(c)$ is homotopic to the image of a cusp by $f_i$.
Since both $k_i$ and $f_i$ are homotopic to $\phi_i$, this means that $c$
is homotopic to a cusp of $S$.
This is a contradiction.

Thus $S'$ must be equal to $S$, and we have a limit pleated surface
$k_\infty \co  S \rightarrow {\bf H}^3/\Gamma_\infty$ touching
$f_\infty(S_0)$.
By Lemma \ref{far surface}, we see that $k_\infty(S) \cap ({\bf
H}^3/\Gamma_\infty)_0$ is contained in $q(E_\mu)$.
Therefore $k_\infty$ can be lifted to a pleated surface $k' \co  S \rightarrow
{\bf H}^3/G$ whose intersection with $({\bf H}^3/G)_0$ is contained in
$E_\mu$.

Now since $k_i$ is homotopic to the pull-back of $k_\infty=q\circ k'$ by an
approximate isometry for sufficiently large $i$, and $k_i$ is homotopic to
$\phi_i$, we see that $k'$ must be homotopic to $\psi$.
As $k_i$ realizes $\overline{\mu}_i$ and $\{\overline{\mu}_i\}$ converges
to $\overline{\mu}$, the pleated surfaces $k_\infty$ and $k'$ realize
$\overline{\mu}$.
As $|\overline{\mu}|=|\mu|$, by changing the transverse measure, $\mu$ can
also be realized by $k'$.
On the other hand, by Lemma~\ref{ending}, $\psi(\mu)$ is an ending
lamination hence $\mu$ cannot be realized by a pleated surface homotopic to
$\psi $.
This is a contradiction.
Thus we have completed the proof of Theorem~\ref{main}.
\endproof

\section{Strong convergence of surface groups}

In Theorems~9.2, 9.6.1 in  Thurston~\cite{Th}, it is stated and roughly
proved that if a sequence of Kleinian groups, which are  isomorphic to a
freely indecomposable Kleinian group (ie satisfying the condition ($*$)
introduced by Bonahon) without accidental parabolics preserving the
parabolicity, converges algebraically to a Kleinian group without
accidental parabolic elements, then the convergence is strong (ie the
geometric limit coincides with the algebraic limit.)
(See also Canary~\cite{Ca}.)
We gave  its detailed proof in Ohshika~\cite{Oh3} except for the case when
the Kleinian group is algebraically isomorphic to a surface group.
The reason why we did not include the case of surface group there is that
it would necessitate to prove that for a convergent sequence, the
hyperbolic structures on the two boundary components cannot degenerate to
the same point in the Thurston boundary.
As this is proved in Theorem~\ref{main}, we can give the proof for the case
of surface group here.
Let $(\Gamma_i,\phi_i)$ be a Kleinian group without accidental parabolic
elements with isomorphism $\phi_i \co  \pi_1(S) \rightarrow \Gamma_i \subset
PSL_2{\bf C}$ for a hyperbolic surface of finite area $S$.
Thurston's original proof in \cite{Th} in this case consists of proving
that the projectivized bending laminations of two boundary components of
the convex cores of ${\bf H}^3/\Gamma_i$ cannot converge to projective
lamination with the same support.
This is exactly the argument on which our proofs of the main theorems are
based. 

\begin{cor}
\label{strong}
Let $S$ be a hyperbolic surface of finite area.
Let $(\Gamma_i,\phi_i)$ be a Kleinian group without accidental parabolic
elements with isomorphism $\phi_i \co  \pi_1(S) \rightarrow \Gamma_i \subset
PSL_2{\bf C}$.
Suppose that $\{(\Gamma_i,\phi_i)\}$ converges algebraically to a Kleinian
group $(G,\psi)$ without accidental parabolic elements.
Then $G$ is also the geometric limit of $\{\Gamma_i\}$.
In other words, $\{\Gamma_i\}$ converges strongly to $G$.
\end{cor}
  
\noindent
\proof
We have only to prove that every subsequence of $\{(\Gamma_i,\phi_i)\}$ has
a subsequence which converges strongly to $(G,\psi)$.
Since a subsequence of $\{(\Gamma_i,\phi_i)\}$ satisfies the condition of
Corollary~\ref{strong}, we only need to show that $\{(\Gamma_i,\phi_i)\}$
in the statement of the corollary has a subsequence strongly converging to
$(G,\psi)$.

By taking a subsequence, we can assume that all of the
$\{(\Gamma_i,\phi_i)\}$ are either quasi-Fuchsian or totally degenerate
groups or totally doubly degenerate groups, and that $\{\Gamma_i\}$
converges geometrically to a Kleinian group $\Gamma_\infty$.
Suppose first that all of the $\{\Gamma_i\}$ are quasi-Fuchsian.
Let $(m_i,n_i) \in {\cal T}(S) \times {\cal T}(S)$ be $Q((\Gamma_i,\phi_i))$.
If both $\{m_i\}$ and $\{n_i\}$ converge inside the Teichm\"{u}ller space
(after taking a subsequence), $\{(\Gamma_i,\phi_i)\}$ converges to a
quasi-Fuchsian group strongly as is well known.
(See for example J{\o}rgensen--Marden~\cite{JM}.)
Assume that one of $\{m_i\}$ and $\{n_i\}$, say $\{m_i\}$, does not
converge inside the Teichm\"{u}ller space and converges to a projective
lamination $[\lambda]$ in the Thurston compactification of the
Teichm\"{u}ller space, and that the other, say $\{n_i\}$, converges inside
the Teichm\"{u}ller space.
Then $G$ is a b--group.
By the same argument as the proof of Lemma~\ref{ending}, the measured
lamination $\lambda$ cannot be realized by a pleated surface homotopic to
$\psi$.
If $\lambda$ is not maximal and connected, as is shown in
Thurston~\cite{Th} or Lemma~4.4 in \cite{Oh1}, $G$ has an accidental
parabolic element, which contradicts our assumption.
Hence $\lambda$ is maximal and connected,  $\psi(\lambda)$ represents an
ending lamination of the geometrically infinite tame end of $({\bf
H}^3/G)_0$, and $G$ is a totally degenerate b--group.
Let $\Sigma_i$ be a boundary component of the convex core of ${\bf
H}^3/\Gamma_i$ corresponding to the ideal boundary component with the
structure $n_i$.
Then as is shown in \cite{Oh3}, the pleated surface $\Sigma_i$ converges
geometrically to a boundary component $\Sigma_\infty$ of the convex core of
${\bf H}^3/\Gamma_\infty$ which can be lifted to a boundary component
$\Sigma$ of the convex core of ${\bf H}^3/G$, which must be the whole
boundary of the convex core as $G$ is a totally degenerate b--group.
Hence a neighbourhood of the geometrically finite end of $({\bf H}^3/G)_0$
is mapped homeomorphically to that of a geometrically finite end of $({\bf
H}^3/\Gamma_\infty)_0$ by the covering projection $q\co  {\bf H}^3/G
\rightarrow {\bf H}^3/\Gamma_\infty$.
On the other hand, by Proposition \ref{covering}, there is also a
neighbourhood of the geometrically tame end of $({\bf H}^3/G)_0$ which is
mapped homeomorphically to a neighbourhood of a  geometrically infinite
tame end of $({\bf H}^3/\Gamma_\infty)_0$ by $q$.
This implies that $G= \Gamma_\infty$.

Next assume that neither $\{m_i\}$ nor $\{n_i\}$ converges inside the
Teichm\"{u}ller space.
After taking a subsequence, we can assume that $\{m_i\}$ converges to a
projective lamination $[\lambda] \in {\cal PL}(S)$ and $\{n_i\}$ converges
to a projective lamination $[\mu] \in {\cal PL}(S)$ in the Thurston
compactification of the Teichm\"{u}ller space.
Since neither $\lambda$ nor $\mu$ can be realized by a pleated surface
homotopic to $\psi$ by Lemma~\ref{ending}, they must be maximal and
connected again by Thurston~\cite{Th} or Lemma~4.4 in \cite{Oh1}.
Then we can apply Theorem~\ref{main} to our situation and see that the
support of $\lambda$ is different from that of $\mu$.
This implies that the end of $({\bf H}^3/G)_0$ with ending lamination
represented by $\psi(\mu)$ is different from one with ending lamination
represented by $\psi(\lambda)$ by Lemma~\ref{unique ending}, hence $G$ is
totally doubly degenerate.
Let $e_\lambda$ and $e_\mu$ denote the two distinct ends of $({\bf
H}^3/G)_0$ with ending laminations represented by $\psi(\lambda)$ and
$\psi(\mu)$ respectively.
By Proposition \ref{covering}, there are neighbourhoods $E_\lambda$ of
$e_\lambda$ and $E_\mu$ of $e_\mu$ such that $q|E_\lambda$ and $q|E_\mu$
are homeomorphisms to neighbourhoods of ends of $({\bf
H}^3/\Gamma_\infty)_0$.
As $({\bf H}^3/G)_0$ has only two ends, this can happen only when
$\Gamma_\infty = G$ or $G$ is a subgroup of $\Gamma_\infty$ of index $2$.
We can see that the latter cannot happen by Lemma~2.3 in \cite{Oh3} (this
fact is originally due to Thurston~\cite{Th}).
Thus we have proved our corollary when all of $\{\Gamma_i\}$ are quasi-Fuchsian.

Next assume that all  the $\Gamma_i$ are totally degenerate b--groups.
Let $m_i$ be the marked hyperbolic structure on $S$ determined by the
conformal structure of $\Omega_{\Gamma_i}/\Gamma_i$, and let $\lambda_i$ be
an  ending lamination of unit length of the geometrically infinite tame end
of $({\bf H}^3/\Gamma_i)_0$.
We can assume that $\{\lambda_i\}$ converges to a measured lamination
$\lambda$ after taking a subsequence.
By the same argument as before, $\lambda$ is maximal and connected, and
$\psi(\lambda)$ represents an ending lamination of $({\bf H}^3/G)_0$ by
Lemma~\ref{ending}.
First assume that $\{m_i\}$ converges inside the Teichm\"{u}ller space.
Then as before, the boundary $\Sigma_i$ of the convex core of ${\bf
H}^3/\Gamma_i$ converges geometrically to a boundary component
$\Sigma_\infty$ of the convex core of ${\bf H}^3/\Gamma_\infty$ which can
be lifted to a boundary component $\Sigma$ of the convex core of ${\bf
H}^3/G$.
Hence $G$ is a totally degenerate b--group, and a neighbourhood of the
geometrically finite end of $({\bf H}^3/G)_0$ is mapped homeomorphically by
$q$ to a neighbourhood of an end of ${\bf H}^3/\Gamma_\infty$.
Then as before, using Proposition \ref{covering}, we can conclude that $G=
\Gamma_\infty$.

Next assume that $\{m_i\}$ does not converge inside the Teichm\"{u}ller space.
Then after taking a subsequence, we can assume that $\{m_i\}$ converges to
a projective lamination $[\mu]$.
By the same argument as before, we can see that $\mu$ is maximal and
connected, and $\psi(\mu)$ represents an ending lamination.
Then by Theorem~\ref{main}, we can see that the support of $\lambda$ is
different from that of $\mu$.
Hence $G$ is totally doubly degenerate, and by Proposition \ref{covering},
we conclude that $G= \Gamma_\infty$.

Finally suppose that all the  $\Gamma_i$ are totally doubly degenerate.
Let $\lambda_i $ and $\mu_i$ be measured laminations of the unit length
such that $\phi_i(\lambda_i)$ and $\phi_i(\mu_i)$ represent ending
laminations of the two geometrically infinite tame ends of $({\bf
H}^3/\Gamma_i)_0$.
By taking a subsequence, we can assume that $\{\lambda_i\}$ converges to a
measured lamination $\lambda$ and $\{\mu_i\}$ converges to a measured
lamination $\mu$ in ${\cal ML}(S)$.
Then as before, both $\lambda$ and $\mu$ are maximal and connected, and
$\psi(\lambda)$ and $\psi(\mu)$ represent ending laminations of $({\bf
H}^3/G)_0$.
By Theorem~\ref{main}, we can see that the support of $\lambda$ is
different from that of $\mu$.
Hence the end of $({\bf H}^3/G)_0$ with ending lamination represented by
$\psi(\lambda)$ is different from that with ending lamination represented
by $\psi(\mu)$ by Lemma~\ref{unique ending}, which implies that $G$ is
totally doubly degenerate.
Then by Proposition \ref{covering} again, we conclude that
$G=\Gamma_\infty$, and the proof is completed.
\endproof

\Addresses\recd
\end{document}